\documentclass[12pt]{article}
\usepackage{amssymb,amsmath, enumitem, mathrsfs}
\usepackage[english]{babel}
\usepackage[utf8]{inputenc}
\usepackage{amsbsy}
\usepackage{rotating}
\usepackage{amsthm}
\usepackage{multirow}
\usepackage{wasysym}
\usepackage{mathrsfs}
\usepackage{verbatim}
\usepackage[colorlinks]{hyperref}
\usepackage{fullpage}
\usepackage{mathdots}
\usepackage{graphicx,subfigure}
\usepackage{stackengine}
\usepackage{tikz,tikz-cd}
\usepackage{tkz-graph}
\usepackage{booktabs}
\usepackage{mathtools}
\usepackage{authblk}
\usepackage{todonotes}
\usetikzlibrary{arrows,positioning,automata}
\PassOptionsToPackage{usernames,dvipsnames,svgnames}{xcolor}
\usepackage{pdflscape}
\usepackage{array}
\newcolumntype{L}{>{$}l<{$}}

\theoremstyle{definition}

\newtheorem{ndefn}{Definition}[section]

\newtheorem{nexap}[ndefn]{Example}

\newtheorem{nalgom}[ndefn]{Algorithm}

\theoremstyle{plain}
\newtheorem*{thm}{Theorem}
\newtheorem{nthm}[ndefn]{Theorem}
\newtheorem*{prop}{Proposition}
\newtheorem{nprop}[ndefn]{Proposition}

\newtheorem{nlem}[ndefn]{Lemma}

\newtheorem{ncor}[ndefn]{Corollary}

\theoremstyle{remark}

\newtheorem{nrk}[ndefn]{Remark}

\DeclarePairedDelimiter\floor{\lfloor}{\rfloor}
\DeclarePairedDelimiter\ab{\lvert}{\rvert}%
\DeclarePairedDelimiter\nm{\lVert}{\rVert}%

\makeatletter
\let\oldabs\ab
\def\ab{\@ifstar{\oldabs}{\oldabs*}}
\let\oldnorm\nm
\def\nm{\@ifstar{\oldnorm}{\oldnorm*}}
\makeatother

\DeclareMathOperator{\cint}{int}
\DeclareMathOperator{\spn}{span}

\DeclareMathOperator{\diag}{diag}

\DeclareMathOperator{\OGr}{OGr}
\DeclareMathOperator{\Flip}{Flip}

\numberwithin{equation}{section}
\numberwithin{figure}{section}

\newcommand{\N}{\mathbb{N}}

\newcommand{\C}{\mathbb{C}}

\newcommand{\Z}{\mathbb{Z}}

\newcommand{\R}{\mathbb{R}}

\newcommand{\x}{\times}

\newcommand{\la}{\langle}
\newcommand{\ra}{\rangle}

\newcommand{\sus}{\subseteq}

\newcommand{\subs}{\subset}

\newcommand{\tl}{\tilde}
\newcommand{\ol}{\overline}

\newcommand{\del}{\Delta}
\newcommand{\lra}{\longrightarrow}
\newcommand{\hra}{\hookrightarrow}

\newcommand{\lmt}{\longmapsto}

\newcommand{\bs}{\backslash}

\newcommand{\mc}[1]{\mathcal{#1}}
\newcommand{\mf}[1]{\mathfrak{#1}}

\newcommand{\mb}[1]{\mathbf{#1}}

\newcommand{\ts}{\textsuperscript}

\newcommand{\lt}{\left}
\newcommand{\rt}{\right}

\newdimen\shadedBaseline\shadedBaseline=-4mm
\newcount\pyramidRow%\newcount\pyramidCol
\newcommand\pyramid[2][\relax]{%
	\ifx\relax#1\relax%
	\else % shade the boxes in #1
	\foreach\bx in {#1} { \filldraw[brown!50!black]\bx+(-.5,-.5) circle (8pt); }
	\fi
	\pyramidRow=0
	\foreach \Row in {#2,...,1} {
		\draw (0,\the\pyramidRow) +(-\Row,0) grid ++(\Row,1);
		\global\advance\pyramidRow by 1  
	}
}

\def\blue#1{{\color{blue} #1}}

\begin{document}

	\title{$DIII$ clan combinatorics for the orthogonal Grassmannian}

	\author[1]{Aram Bingham}
	\author[2]{\"Ozlem U\u{g}urlu}
	
	\affil[1]{{\small Tulane University, New Orleans; abingham@tulane.edu}}    
	\affil[2]{{\small Palm Beach State College, Boca Raton; ugurluo@palmbeachstate.edu}}
	
	\normalsize

	\date{\today}
	\maketitle
	
	\begin{abstract}
		
Borel subgroup orbits of the classical symmetric space $SO_{2n}/GL_n$ are parametrized by $DIII$ $(n,n)$-clans. We group the clans into ``sects'' corresponding to Schubert cells of the orthogonal Grassmannian, thus providing a cell decomposition for $SO_{2n}/GL_n$. We also compute a recurrence for the rank polynomial of the weak order poset on $DIII$ clans, and then describe explicit bijections between such clans, diagonally symmetric rook placements, certain pairs of minimally intersecting set partitions, and a class of weighted Delannoy paths.  Clans of the largest sect are in bijection with fixed-point-free partial involutions.
		\vspace{.2cm}
		
		
	\end{abstract}

\section{Introduction}\label{S:Introduction}

% are two main approaches to the cohomology of a flag variety $G/B$ of a complex reductive algebraic group $G$ with Borel subgroup $B$. The first is the Borel picture \cite{borel53}, which identifies $H^*(G/B)$ with the the coinvariant algebra of a reflection representation of the Weyl group of $G$.\footnote{See also \cite[Chapter 4]{hiller} for a more recent account in the complex algebraic group setting.} The other, perhaps more in the spirit of the original Schubert calculus, proceeds from the Bruhat decomposition of $G$ which, when interpreted in terms of $B$-orbits in the flag variety, provides a cell decomposition of $G/B$. The closures of these cells are called \emph{Schubert varieties}, and their Poincar\'e duals give a basis for the integral cohomology of $G/B$. This paper is a continuation of a program, outlined in \cite{aoSects}, to take the latter approach toward studying the cohomology of certain \emph{symmetric spaces} (or \emph{symmetric varieties}).

Symmetric spaces are an important class of \emph{spherical} $G$-varieties. If $G$ is a complex reductive algebraic group, spherical varieties are those for which a Borel subgroup $B$ acts with finitely many orbits. The theory of spherical varieties encompasses that of toric varieties, and their classification can be given similarly in terms of ``colored fans'' (see \cite{knop} for an introduction). As such, these varieties present rich opportunities for combinatorial investigation as complement to their significance in algebraic geometry and representation theory. 

We define symmetric spaces as follows. Supposing $G$ has an algebraic automorphism $\theta$ of order two, then the fixed-point subgroup $L:=G^\theta$ is called a \emph{symmetric subgroup} and the quotient $G/L$ is the associated symmetric space.\footnote{In the literature, this definition is often broadened to include any space $G/K$ where $K$ lies between the connected component of the identity of $G^\theta$ and its normalizer. } For simple $G$, symmetric spaces were effectively classified by Cartan in the course of classifying real forms of simple Lie algebras over the complex numbers (see \cite[Chapter 10]{helgason}).

%Further, the bulk of the theory works over any algebraically closed field of characteristic other than 2 \cite{richardsonSpringer1990}.

Within the classification, there are a few cases which are closely related to \emph{Grassmannian varieties}, which are realized as homogeneous spaces $G/P$ where $P$ is a maximal parabolic subgroup of $G$. Grassmannian varieties parameterize vector subspaces of a given vector space, and their cohomology theory is  both a classical subject \cite{kleiman} and a central topic in modern algebraic combinatorics. For the symmetric spaces in question, the subgroup $P$ admits a \emph{Levi decomposition} as $P=L\ltimes R_u(P)$ where $L=G^\theta$ as before, and $R_u(P)$ is the unipotent radical of $P$ \cite[Section 30.3]{humphreysLAG}.

This paper concerns the third of three cases in which this occurs,\footnote{See Remark~\ref{rk:BDI}.} namely the symmetric space of type \emph{DIII}. Similar analysis was performed for type $AIII$ ($SL_{p+q}/S(GL_p \x GL_q)$) in~\cite{bcSects} and \cite{canGenesis}, and for type $CI$ ($Sp_{2n}/GL_n$) in \cite{aoSects}. The labels come from Cartan's original classification, which can be viewed as a refinement of the classification of simple Lie algebras over the complex numbers. Type $DIII$ refers to the quotient $SO_{2n}/GL_n$. A realization in coordinates will be given in the next section. Note that while all of the symmetric space theory we use is valid over an algebraically closed field of characteristic other than two \cite{richardsonSpringer1993}, all groups in this paper are taken to be over the complex numbers.

Borel orbits in symmetric spaces are often parameterized by sets of \emph{clans}, following terminology of Matsuki-Oshima \cite{matsukiOshima}. Since the work of Yamamoto \cite{yamamoto}, clans often appear as strings of $+$ and $-$ symbols interspersed with pairs of matching natural numbers, for example ${+}12{+}{-}12{-}$ (see Definition~\ref{def:clans}). As $B$-orbits in $G/L$ are in bijection with $L$-orbits in $G/B$, a clan encodes the data of a representative flag for the corresponding $L$-orbit in the flag variety $G/B$, but they may also be regarded as \emph{signed involutions} of the symmetric group (see Definition~\ref{def:signedinv}). This paper will describe some of the geometry of $SO_{2n}/GL_n$ in terms of clans, as well as provide some relevant combinatorial results.

We shall now describe the organization of this paper. From now on, let $\theta$ be an involution on $G:=SO_{2n}$ which has $L:=G^\theta\cong GL_n$ as fixed-point subgroup, and let $B$ be a Borel subgroup of $G$ containing a $\theta$-stable maximal torus $T$ of $G$. We will refer to the clans which parametrize the $B$-orbits in $SO_{2n}/GL_n$ as \emph{$DIII$ $(n,n)$-clans} (see Definition~\ref{defn:Dclans}), or just \emph{$DIII$ clans} if $n$ is either clear from context or irrelevant. After setting down some notation and terminology in Section~\ref{S:Notation}, our first result is Theorem~\ref{thm:Dflag} which provides flags to represent $GL_n$-orbits in $SO_{2n}/B$, using results of \cite{wyserThesis}. These particular flags had been overlooked in previous literature on clans.

The symmetric subgroup $GL_n \subs SO_{2n}$ can be realized as the Levi factor of a maximal parabolic subgroup $P$ such that $SO_{2n}/P$ is $\OGr(n, \C^{2n})$, the \emph{orthogonal Grassmannian} of maximal ($n$-dimensional) isotropic subspaces of $\C^{2n}$. This gives us a canonical $G$-equivariant projection map $\pi: SO_{2n}/GL_n \to \OGr(n, \C^{2n})$. Borel orbits in Grassmannian varieties are called \emph{Schubert cells}. A \emph{sect} is a collection of clans indexing $B$-orbits which map to the same Schubert cell under $\pi$. In Theorem~\ref{decomp}, we prove a description of the sects of $SO_{2n}/GL_n$ which matches that previously given by the authors and Can for types $AIII$ and $CI$. 

From results of~\cite{bcSects}, the sects provide a cell decomposition and $\Z$-basis for the (co)homology of $SO_{2n}/GL_n$. Further, an isomorphism of cohomology rings $H^*(G/L)\cong H^*(G/P)$ follows from the fact that the fibration
\[\C^a \cong R_u(P) \cong P/L \hra G/L \lra G/P \]
gives $\pi:G/P \to G/L$ the structure of an affine bundle, and applying the Leray-Hirsch theorem. The latter ring is understood to be the subring of $W_L$-invariants within coinvariant algebra of a reflection representation of the Weyl group of $SO_{2n}$, where $W_L$ is the Weyl group of $L$.\footnote{See \cite[Chapter 4]{hiller} for background; note $W_L$ is also the Weyl group of $P$.} 

Clans form a graded poset under the \emph{weak order}, first defined in \cite{richardsonSpringer1990}. The covering relations of this poset are given by the action of minimal standard parabolic subgroups $P_{s}$ on corresponding $B$-orbits, where $s$ is a simple reflection of the Weyl group $W=N_G(T)/T$. We recall a combinatorial description of this order and its associated \emph{length function} to compute the following recurrence relation for the rank polynomial of this weak order poset in Section~\ref{S:PartialOrder}.
\begin{thm}[\ref{thm:lgf}]
The rank polynomial $A_n(t)$ of the weak order poset on $DIII$ clans satisfies the following recurrence relation:
\begin{equation*} 
A_n(t)=2A_{n-1}(t) + (t+ t^2 + \cdots t^{n-2}+2t^{n-1}+t^n +\cdots +t^{2n-3} )A_{n-2}(t).
\end{equation*}
\end{thm}

%$DIII$ $(n,n)$-clans are a subset of the set of all $(n,n)$-clans (which parameterize Borel orbits in the type $AIII$ case). It turns out that the additional restrictions on $DIII$ clans give combinatorial coincidences with other well-studied families of objects, and
This recurrence easily gives a generating function and recurrence for $\del_n$, the number of $DIII$ $(n,n)$-clans, but we also obtain an explicit formula by a different method in Section~\ref{S:Formula}.
\begin{prop}[\ref{cor:Dformula}]
	The number of $DIII$ $(n,n)$-clans  is 
	\begin{equation*}
%	\label{eq:Deln}
	\del_n= \sum_{r=0}^{\floor{\frac{n}{2}}} 2^{n-2r-1}\frac{n!}{r!(n-2r)!}.
	\end{equation*}
\end{prop}
The rest of Section~\ref{S:Formula} describes bijections between $DIII$ clans and other combinatoral families of objects. The first involves the number of inequivalent placements of $2n$ non-attacking rooks on a $2n \x 2n$ board with symmetry across each of the main diagonals, which was written about in the classic text of Lucas~\cite{lucas}. A bijection between $DIII$ $(n,n)$-clans and such rook placements is given in Section~\ref{subsec:rookpyr}, by extracting a triangular portion of the square board and analyzing this resulting ``pyramid.''  These pyramids also make it easy to describe a (near) bijection with objects studied by Pittel in~\cite{pittel}. Precisely, these are ordered pairs $(p,p')$ of partitions of an $n$-element set such that $p$ consists of exactly two blocks. This map is described in Section~\ref{subsec:setparts}. 

 Schubert cells of $G/P$ can be parameterized by lattice paths, which are also a tool for understanding their geometry. As a step towards extending these ideas to the symmetric space above, we present a bijection between $DIII$ $(n,n)$-clans and certain \emph{weighted Delannoy paths} in Section~\ref{subsec:LP}. We do not further investigate the classes of $B$-orbit closures in the cohomology of $G/L$, but it is our hope that a lattice path model for the orbits may be helpful for the future development of tableaux combinatorics to describe multiplication in the cohomology ring of $G/L$, extending the Littlewood-Richardson rule. Wyser has related clan orbit closures to Richardson subvarieties of flag varieties in order to extract some information on Schubert calculus of flag varieties \cite{wyserCD, wyserRich}.
 
Finally, we look at the pre-image of the dense Schubert cell of $SO_{2n}/GL_n$, which we call the \emph{big sect}. We prove that the clans of the big sect are in bijection with the set of partial fixed-point-free involutions of an $n$-element set, denoted by $\mc{PF}_n$. The elements of $\mc{PF}_n$ parametrize congruence orbits of the upper triangular invertible matrices on the set of skew-symmetric matrices, as described in~\cite{cherniavskySkew}. Equipped with the closure order of the orbits of that action, they form a poset which has also been studied in~\cite{cantim}. Proof of the following theorem will appear in the first author's Ph.D. thesis. 
\begin{thm}
The closure order on $DIII$ (n,n)-clans of the largest sect is isomorphic to the poset of partial fixed-point-free involutions on $n$ letters with the congruence orbit closure order.
\end{thm}

\begin{nrk}\label{rk:BDI}
The list of symmetric spaces with symmetric subgroup equal to a Levi subgroup often includes the \emph{type $BDI$} spaces $SO_n/(SO_2\x SO_{n-2})$, rounding out the symmetric spaces of \emph{Hermitian type}. However, our definition of symmetric space (which matches that of \cite{richardsonSpringer1990,richardsonSpringer1993}) excludes this case from consideration. Some details on Hermitian-type spaces (including type $BDI$) using an alternative orbit parametrization can be found in \cite{richardsonSpringer1993, ryan}.
\end{nrk}

\section{Notation and preliminaries}\label{S:Notation}

Let $n$ be a positive integer.
First, we describe our realization of $SO_{2n}$. We follow most of the notations of~\cite{wyserThesis}. 

Let $J_{2n}$ denote the $2n \x 2n$ matrix with 1's along the anti-diagonal and 0's elsewhere. Then, we set 
\[ G:=SO_{2n}= \{ g\in SL_{2n} \mid g^t J_{2n} g = J_{2n}  \}. \]
Let $\cint(g): SL_{2n} \to SL_{2n}$ denote the map defined by 
\[ \cint (g)(h)= ghg^{-1} .\]
Now define the matrix 
\[I_{n,n} := \begin{pmatrix}
I_n & 0 \\
0 & -I_n 
\end{pmatrix}. \]
Then we check that we have an involution $\theta$ on $SL_{2n}$ defined by $\theta:=\cint(iI_{n,n}).$

Since $(iI_{n,n})^{-1}=-i I_{n,n}$, then if 
\[g=\begin{pmatrix}
A & B\\ 
C & D
\end{pmatrix}\]
is the $n \x n $ block form of $g$, we have

\[\theta(g)=
\begin{pmatrix}
i I_n & 0 \\ 
0 & -i I_n
\end{pmatrix}
\begin{pmatrix}
A & B \\
C & D 
\end{pmatrix}
\begin{pmatrix}
-i I_n & 0 \\ 
0 & i I_n
\end{pmatrix}
= \begin{pmatrix}
A & -B \\
-C & D
\end{pmatrix}. \]
Observe that the restriction of $\theta$ to $SO_{2n}$ induces an involution on that group as well. The fixed points of this involution must be block diagonal, that is
\[ \theta(g)=g \iff g= \begin{pmatrix}
A & 0 \\
0 & D
\end{pmatrix}.\]
Furthermore, membership in the special orthogonal group forces $D=J_n(A^{-1})^tJ_n$. Thus, $A$ can be any invertible $n \x n$ matrix and this completely determines $g$, so the fixed point subgroup $L$ is isomorphic to $GL_n$.

Next, we fix some combinatorial notation. We will write $\mc{S}_n$ for the symmetric group of permutations of $[n]:=\{1,\dots, n\}$. If $\pi \in \mc{S}_n$, then its one-line notation is the string $\pi_1\pi_2\dots \pi_n$, where $\pi_i = \pi(i)$ for $1\leq i \leq n$. For instance, $\pi=164578329$ is the one-line notation for the permutation $\pi \in \mc{S}_9$ with cycle decomposition $(1) (2\ 6\ 8)(3\ 4\ 5\ 7) (9)$.

An involution is an element of $\mc{S}_n$ of order at most two, and the set of involutions in $\mc{S}_n$ is denoted by $\mc{I}_n$. An involution $\pi\in\mc{I}_n$ can be written in cycle notation in the canonical form
$$
\pi = (a_1\ b_1)(a_2\ b_2) \dots (a_k\ b_k) (d_1)\dots (d_{n-2k}),
$$ 
where $a_i < b_i$ for all $1 \leq i \leq k$, $a_1 < a_2 < \dots < a_k$, and $d_1<\dots < d_{n-2k}$. \emph{Signed involutions} are involutions where the fixed points are decorated with a choice of sign, $+$ or $-$.
\begin{ndefn} \label{def:signedinv}
	Let $p$ and $q$ be positive integers where $q \leq p$. A \emph{signed $(p,q)$-involution} is a signed involution on ${p+q}$ letters such that there are $p - q $ more $+$'s than $-$'s. 
\end{ndefn}
For example, $\pi = (1\; 8)(2\; 4)(3^+)(5^-)(6^-)(7^+)$ is a signed $(4,4)$-involution. Clans can be thought of as an alternative presentation of signed involutions.
\begin{ndefn} \label{def:clans}
	Let $p$ and $q$ be two positive integers where $q\leq p$. A $(p,q)$\emph{-clan} $\gamma$ is a string of $p+q$ symbols from $\N\cup\{+,-\}$ such that 
	\begin{enumerate}
		\item there are $p-q$ more $+$'s than $-$'s;
		\item if a natural number appears in $\gamma$, 
		then it appears exactly twice. 
	\end{enumerate} 
	For example, $12{+}2 1$ is a $(3,2)$-clan and ${+}1{+}1$ is a $(3,1)$-clan. Clans $\gamma$ and $\gamma'$ are considered to be equivalent if the positions of the matching number pairs are the same in both clans. For example, $\gamma=1 1 22$ and $\gamma'=22 1 1$ are the same $(2,2)$-clan, since both of $\gamma$ and $\gamma'$ have matching numbers in positions $(1,2)$ and in positions $(3,4)$. 
\end{ndefn}
Evidently, a clan is just a signed involution in list notation where 2-cycles give the positions of matching natural numbers and the locations of fixed points are occupied by their signatures. To illustrate the equivalence between signed involutions and clans, we remark that the signed $(4,4)$-involution $(1\; 8)(2\; 4)(3^+)(5^-)(6^-)(7^+)$ can be regarded as the $(4, 4)$-clan $12{+}2{-}{-}{+}1$.

Throughout this paper we prefer to use clans, though we will occasionally like to refer to the \emph{underlying involution} of a clan, which is obtained by simply ignoring the signs on fixed-points in the associated signed involution. We will denote the underlying involution of clan $\gamma$ by $\sigma_\gamma$. In a clan $\gamma$, the matching natural numbers of a pair coming from a transposition in $\sigma_\gamma$ will be referred to as \emph{mates} of one another.  %Also observe that $p$ is equal to the number of fixed points in $\pi$ with a $+$ sign attached plus the number of two-cycles in $\pi$, while $q$ is equal to the number of fixed points in $\pi$ with a $-$ sign attached plus the number of two-cycles in $\pi$. 

Next we will specify the $DIII$ clans. Let $\gamma$ be a 
clan of the form $\gamma=c_1\cdots c_n$.
The \emph{reverse} of $\gamma$,
denoted by $rev(\gamma)$,
is the clan 
$$
rev(\gamma) := c_nc_{n-1}\cdots c_1.
$$
We obtain the \emph{negative} of $\gamma$, denoted by $\ol{\gamma}$, by changing all $+$'s in $\gamma$ to $-$'s, and vice versa, leaving the natural numbers unchanged.
Now, we define symmetric and skew-symmetric clans. 
\begin{ndefn}
	A $(p,q)$-clan $\gamma$ is called \emph{symmetric} if 
	$$\gamma =  rev(\gamma),$$
	and is called \emph{skew-symmetric} if 
	$$\gamma =   {rev(\ol{\gamma})}.$$
\end{ndefn} 
\begin{nexap}
	Consider the clan $\gamma = {+}{-}1 2 3 3 1 2{+}{-}$. Its reverse is $ {-} {+} 2 1 3 3 2 1 {-} {+}$. Since $\gamma = {rev (\ol\gamma)}$, it is a skew-symmetric $(5, 5)$-clan.
	
	The clan $\tau = 1 2 3 4 54 5 3 2 1$ is a skew-symmetric $(5,5)$-clan which is also symmetric, as it contains no $\pm$ symbols. 
\end{nexap}
A pair of mates $(c_i,c_j)$ exchanges places with another pair of mates $(c_{2n+1-j},c_{2n+1-i})$ upon taking the reverse of a clan. Such pairs shall be called \emph{opposing pairs}. For instance $(c_3, c_7)$ and $(c_4,c_8)$ are opposing pairs in $\gamma = {+}{-}1 2 3 3 1 2{+}{-}$.
\begin{ndefn} \label{defn:Dclans}
The set of \emph{$DIII$ $(n,n)$-clans} consists of those $(n,n)$-clans $\gamma=c_1\cdots c_{2n}$ which satisfy following the additional conditions:
\begin{enumerate}
	\item $\gamma$ is skew-symmetric, that is $\gamma=rev(\ol{\gamma})$; 
	\item if $c_i\in \N$, then $c_i \neq c_{2n+1-i}$;
	\item the total number of $-$'s and pairs of matching natural numbers among $c_1\cdots c_n$ is even.
\end{enumerate}
\end{ndefn}
Recall that the Coxeter group of type $D_n$, which is the Weyl group of $SO_{2n}$, can be regarded as the signed permutations on $n$ letters with an even number of sign changes. This can also be viewed as a subgroup of $S_{2n}$ by identifying the symbol $-i$ with $2n+1-i$ for each $1\leq i \leq n$. The underlying involutions of $DIII$ clans are then involutions of a type $D_n$ Weyl group. Note that unlike $(n,n)$-clans for the type $AIII$ and $CI$ symmetric spaces, not all involutions of the Weyl group are attainable as the underlying involution of some clan. In particular, condition $2$ of Definition~\ref{defn:Dclans} prohibits the longest element of the type $D_n$ Weyl group from arising as an underlying involution.\footnote{When $n$ is even, the longest element takes $i\mapsto -i$ for all $1\leq i \leq n$, which would be underlying the clan $12\cdots n n \cdots 2 1$.}

We shall write $\del(n)$ to denote the set of $DIII$ $(n,n)$-clans. It was first stated in~\cite{matsukiOshima} and proved in~\cite{wyserThesis} that $DIII$ $(n,n)$-clans parametrize $L$-orbits in $SO_{2n}/B$. Our notation for clans comes from the latter source. In the next section, we will produce representative flags for each orbit and describe sects for these clans.

\section{Sects}\label{S:sects}
\subsection{Background}
In order to describe sects and representative flags for $DIII$ clans, we must first visit the theory of parabolic subgroups of special orthogonal groups. We refer the reader to~\cite{malleTesterman} for more details.

Given a vector space $V$ with bilinear form $\omega$, recall that an \emph{isotropic subspace} $W$ is one such that $\omega(\mb{u},\mb{v}) = 0$ for all vectors $\mb{u}, \mb{v} \in W$. If we also use $\omega$ to stand for the matrix which represents this bilinear form in a particular choice of basis, this condition becomes  $\mb{u}^t \omega \mb{v} =\mathbf{0}$. A \emph{polarization} of $V$ is then a direct sum decomposition of $V$ into subspaces which are each isotropic (with respect to $\omega$), that is
$V = V_- \oplus V_+$.

Let $E_n$ be the subspace generated by standard basis vectors $\{\mb{e}_i \mid 1\leq i \leq n\}$. It is easy to check that this is an isotropic subspace of $\C^{2n}$ with respect to $J_{2n}$, and in fact it is maximal with respect to inclusion of isotropic subspaces. There is a distinguished polarization of $\C^{2n}$ as 
\begin{equation} \label{eqn:polarization}
    V=E_n\oplus \tl{E}_n,
\end{equation}
where $\tl{E}_n$ is the subspace spanned by $\{\mb{e}_{n+1},\dots,\mb{e}_{2n}\}$. Note, however, that (infinitely) many other isotropic subspaces could replace $\tl{E}_n$ in the direct sum decomposition above. 

An \emph{isotropic flag} is defined as a sequence of vector spaces
\[ \{\mathbf{0} \} \subs V_1 \subs \dots V_r \subs V \]
such that $V_i$ is an isotropic subspace of $V$ for all $1\leq i \leq r$.

Taking $V=\C^{2n}$, we have a bilinear form given by the matrix $J_{2n}$ used to define our realization of $G:=SO_{2n}$. From~\cite[ Proposition 12.13]{malleTesterman}, the parabolic subgroups of $SO_{2n}$ are precisely the stabilizers of flags which are isotropic with respect to the form $J_{2n}$.\footnote{Malle and Testerman define their special orthogonal group by a form which is a scalar multiple (by one-half) of the one presented here, but the resulting group that leaves the form invariant is the same.}

The stabilizer of the the flag $\{\mb{0}\} \subs E_n \subs \C^{2n}$ is the parabolic subgroup $P$ consisting of matrices with $n \x n$ block form
\begin{equation}\label{eq:PL}
  \begin{pmatrix}
* & * \\
0 & * 
\end{pmatrix}
\;
\text{and which has Levi factor }
L =\lt\{ \begin{pmatrix}
A & 0 \\
0 & J_n(A^{-1})^t J_n 
\end{pmatrix} \;\middle|\; A \in GL_n \rt\}; 
\end{equation}
see~\cite[p. 144]{malleTesterman} or \cite[Section 8.1]{garrettBuildings} for related discussion. Thus, we see that the Levi subgroup of this parabolic subgroup coincides with a symmetric subgroup of type $DIII$, that is $G^\theta=L$ where $\theta$ is the involution of Section~\ref{S:Notation}. Furthermore, the subgroup $L$ is exactly the stabilizer of the polarization of (\ref{eqn:polarization}). The association of the symmetric pair $(G,L)$ with a polarization is another feature that type $DIII$ has in common with type $AIII$ and $CI$ (see~\cite[p. 511]{goodmanWallach}).

The upshot of this coincidence is that we have a $G$-equivariant projection map 
\[\pi: G/L \lra G/P\]
which we can analyze. Letting $B$ be the Borel subgroup of upper triangular matrices in $G$ \cite[p. 39]{malleTesterman} (which contains the $\theta$-stable maximal torus of diagonal matrices) we can relate the $B$-orbits in $G/P$, which are Schubert cells, to the $B$-orbits in $G/L$. The equivariance of $\pi$ allows us to ask precisely which clans constitute the pre-image of a particular Schubert cell. We call such a collection of clans the \emph{sect} associated to the Schubert cell.

In the literature, clans usually parametrize symmetric subgroup orbits in a flag variety by encoding the information of how flags in that orbit relate to a reference polarization. In type $DIII$, one may consider $L$-orbits in $G/B$, which can be identified with one component of the variety of all full flags isotropic with respect to $J_{2n}$. For $J_{2n}$, a \emph{full} isotropic flag $V_\bullet$ in $\C^{2n}$ is a sequence of vector subspaces $\{V_i\}_{i=0}^n$ such that
\begin{equation}\label{eq:isoflag}
\{\mb{0}\}= V_0 \subs V_1 \subs V_2 \subs \dots \subs V_n, 
\end{equation} 
where $\dim V_i=i$ for all $i$ and $V_n$ is a maximal isotropic subspace. We find it convenient to write 
\[ V_\bullet = \la \mb{v}_1, \mb{v}_2, \dots, \mb{v}_n  \ra \] to indicate that $V_\bullet$ is the flag with $V_i= \spn \{ \mb{v}_1, \dots, \mb{v}_i \}$ for all $1\leq i \leq n $. Any full isotropic flag is canonically extended to a full flag in $\C^{2n}$
\[\{\mb{0}\} \subs V_1 \subs \dots \subs V_{2n-1} \subs V_{2n}=\C^{2n} \]
by assigning 
\[V_{2n-i} := V_i^\perp = \{ \mb{v}\in \C^{2n} \mid \omega(\mb{v},\mb{w})= 0,\  \forall \mb{w}\in V_i \}, \]
so we may abuse notation slightly by using $V_\bullet$ to refer to either one. For instance, the standard isotropic full flag $E_\bullet:= \la \mb{e}_1, \dots, \mb{e}_n \ra$ can be written in extended form as
\[ E_\bullet= \la \mb{e}_1, \dots, \mb{e}_n, \mb{e}_{n+1},\dots, \mb{e}_{2n} \ra.\]

If $g\in G$ is a matrix whose $i$\ts{th} column is a vector $\mb{v}_i$, then the isotropic flag corresponding to the coset $gB$ will be given by $V_\bullet=\la \mb{v}_1, \dots, \mb{v}_{2n}\ra$, and vice versa. For example, the coset of the identity matrix $I_{2n}$ corresponds to the standard isotropic flag $E_\bullet$.

The space of full flags isotropic with respect to $J_{2n}$ is a disconnected double cover of $G/B$; it consists of two isomorphic $SO_{2n}$-orbits. Since we will represent flags by matrices $g$ that identify cosets $gB\in G/B$, and we want the standard flag 
$E_\bullet$ to identify with the identity coset, then to guarantee that a $J_{2n}$-isotropic flag $V_\bullet$ is in the same $SO_{2n}$ orbit as $E_\bullet$ we must add the additional condition that $\dim (V_n \cap E_n)\equiv n \bmod 2$ \cite[p. 106]{wyserThesis}. The set of such flags is then an honest homogeneous space for $SO_{2n}$.

We must present a few more definitions before describing the process of obtaining orbit-representative flags; see also~\cite{bcSects}.
\begin{ndefn}
	Given an $(n,n)$-clan $\gamma= c_1\cdots c_{2n}$, one obtains the \emph{default signed clan} associated to $\gamma$ by assigning to $c_i$ a ``signature'' of $-$ and to $c_j$ a ``signature'' of $+$ whenever $c_i=c_j \in \N$ and $i< j$. We denote this default signed clan by $\tl{\gamma} =\tl{c}_1 \cdots  \tl{c}_{2n}$. 
\end{ndefn}

For instance, $\tl{\gamma} = + 1_- 2_- 1_+ 2_+ -$ is the default signed clan of $\gamma=+ 1 2 1 2 -$. Note that the signature of $\tl{c}_i$ is just the symbol itself in case $c_i$ is $+$ or $-$.
\begin{ndefn} Given a default signed clan $\tl{\gamma}$, define a permutation ${\sigma} \in \mc{S}_{2n}$ which, for $i\leq n$:
	\begin{itemize}
		\item assigns $\sigma(i)=i$  and $\sigma(2n+1-i)= 2n+1-i$ if $\tilde{c}_i$ is a symbol with signature $+$.
		\item assigns $\sigma(i)=2n+1-i$ and $\sigma(2n+1-i)=i$ if $\tl{c}_i$ is a symbol with signature $-$.
	\end{itemize} We call ${\sigma}$ the \emph{default permutation} associated to $\gamma$. 
\end{ndefn}
Returning to our example, $+1212-$ has default permutation $1 5 4 3 2 6$ (in one-line notation). Note that $\sigma$ is an involution, and it is the $\sigma'$ which results from choosing $\sigma''=id$ in the context of~\cite[Theorem 3.2.11]{yamamoto}.

\subsection{Sects for $DIII$ clans} \label{subsec:flagerr}
Fix, as before, $G=SO_{2n}$, $B$ its Borel subgroup of upper triangular matrices, and $P$ and $L$ as defined by (\ref{eq:PL}).  Next, we show how to obtain representative flags for $L$-orbits in $G/B$ from $DIII$ $(n,n)$-clans using a variant of the methods in~\cite{yamamoto}. To ensure that we obtain a complete set of representative flags for all $DIII$ $(n,n)$-clans, we apply the following instance of~\cite[Theorem 1.5.8]{wyserThesis}.

\begin{nthm} \label{thm:DPar}For the symmetric pair $(G,L)=(SO_{2n}, GL_n)$ of type $DIII$, each $L$-orbit of $G/B$ is equal to the intersection of an $S(GL_n \x GL_n)$-orbit in the flag variety $X'$ of $SL_{2n}$ with the isotropic flag variety, viewed as a subvariety $X \sus X'$.
\end{nthm} 
This theorem accords with the view of $DIII$ $(n,n)$-clans as a subset of all $(n,n)$-clans whose elements satisfy extra conditions; the inclusion of sets of clans is reflected in the containment of the respective orbits. Then for each $DIII$ $(n,n)$-clan $\gamma$, to obtain a representative flag for the $L$-orbit $Q_\gamma$, it suffices to produce an isotropic flag $V_\bullet(\gamma)$ which satisfies 
\[\dim (V_n(\gamma) \cap E_n) \equiv n \bmod 2\] and can also be produced by~\cite[Theorem 2.2.14]{yamamoto}, as that theorem provides flags for type $AIII$ clans. This will give us a full set of representative flags on which we can perform the sect analysis.  

\begin{nthm} \label {thm:Dflag}Given a $DIII$ $(n,n)$-clan $\gamma=c_1 \cdots c_{2n}$ with default permutation ${\sigma}$, define a flag $V_\bullet(\gamma)=\la \mb{v}_1, \dots, \mb{v}_{2n}\ra$ by making the following assignments.
	
	\begin{itemize}
		\item If $c_i = \pm $, set \[v_i = \mb{e}_{\sigma(i)}. \]
		\item If $c_i=c_j \in \N$ where $i< j$, so that $c_{2n+1-i}=c_{2n+1-j}$ as well, with $i<2n+1-j$, then set 
		\begin{align*}
		\mb{v}_i &= \frac{1}{\sqrt{2}}(\mb{e}_{\sigma(i)}+\mb{e}_{\sigma(j)}), \\
		\mb{v}_j &=\frac{1}{\sqrt{2}}(\mb{e}_{\sigma(i)}-\mb{e}_{\sigma(j)}), \\
		\mb{v}_{2n+1-i} &= \frac{1}{\sqrt{2}}(\mb{e}_{\sigma(2n+1-i)}+\mb{e}_{\sigma(2n+1-j)}), \\
		\mb{v}_{2n+1-j} &= \frac{1}{\sqrt{2}}(\mb{e}_{\sigma(2n+1-i)}-\mb{e}_{\sigma(2n+1-j)}).
		\end{align*}
	\end{itemize}
	Then $V_\bullet(\gamma)$ is a representative flag for the $L$-orbit $Q_\gamma$. Furthermore, if $g_\gamma$ is the matrix defined by letting $\mb{v}_i$ be its $i$\ts{th} column, then $Q_\gamma=Lg_\gamma B/B$ in $G/B$. Matrices/flags obtained in this way from all $DIII$ $(n,n)$-clans constitute a full set of representative flags for $L$-orbits in $G/B$.
\end{nthm}
\begin{proof} The verification that $g_\gamma \in SO_{2n}$ is routine (if tedious) linear algebra. From this it follows that $V_\bullet(\gamma)$ is isotropic with respect to $J_{2n}$. 
	
	Next we argue that
	\[\dim (V_n(\gamma) \cap E_n) \equiv n \bmod 2.\]
	If a $-$ appears at $c_i$, then $\mb{v}_i=\mb{e}_r$ for some $r>n$. Thus, for each $-$ among the first $n$ symbols, \hbox{$\dim (V_n(\gamma) \cap E_n)$} is reduced by one (compared to when $V_n=E_n$). For each $c_i=c_j \in \N$, the vector subspace spanned by $\mb{v}_i$ and $\mb{v}_j$ is equal to $\spn(\mb{e}_q, \mb{e}_r)$ for some $q\leq n$ and some $r>n$. Then, each pair of matching natural numbers among $c_1,\dots, c_n$ reduces $\dim (V_n(\gamma) \cap E_n)$ by one as well. Since there are an even number of $-$'s and pairs of matching natural numbers among the first $n$ symbols, 
	\[\dim (V_n(\gamma) \cap E_n)= n-2k\equiv n \bmod 2,\]
	for some natural number $k$.
	
	Finally, we mention how to obtain this flag from~\cite[Theorem 2.2.14]{yamamoto}.
	\begin{ndefn}
	Let a \emph{family} in $\gamma=c_1\cdots c_{2n}$ mean a collection of symbols $c_i, c_j, c_{2n+1-j}, c_{2n+1-i}$ with $c_i=c_j\in \N$, $i<j$, and $i< 2n+1-j$.
	\end{ndefn} 
	For each family in $\gamma$, modify the default signed clan by flipping the signatures of $\tl{c}_i$ and $\tl{c}_j$ so that they are $+$ and $-$, respectively. Denote the signed clan so obtained by $\tl{\gamma}'$. To reflect this adjustment, we also modify the default permutation $\sigma$ by swapping $\sigma(i)$ and $\sigma(j)$ for all such families. Denote the permutation so obtained by $\sigma'$. 
	
Then $\sigma'$ satisfies the conditions of~\cite[Theorem 2.2.14]{yamamoto}, and we claim the flag produced by that theorem using $\tl{\gamma}'$ and $\sigma'$, and denoted by $V'_\bullet(\gamma)=\la \mb{v}_1',\dots, \mb{v}_{2n}' \ra$, is the same flag as $V_\bullet(\gamma)$ constructed above. Indeed, applying that theorem, one finds $\mb{v}_i'=\mb{v}_i=\mb{e}_{\sigma(i)}$ whenever $c_i= \pm$, and for any family $c_i,c_j, c_{2n+1-j}, c_{2n+1-i}$, we get
	\[\mb{v}_i=\mb{v}_i' \quad \text{and}\quad \mb{v}_j=\mb{v}_j' \quad \text{and} \quad \mb{v}_{2n+1-j}=-\mb{v}_{2n+1-j}' \quad \text{and} \quad \mb{v}_{2n+1-i}=\mb{v}_{2n+1-i}'.\]
	Clearly, these generate the same flags. Thus, the theorem is proved.
\end{proof}
For example, the matrix representative for the clan $\gamma=+1212- $ from Theorem \ref{thm:Dflag} is
\[ g_\gamma =
\begin{pmatrix}
1 & 0 & 0 & 0 & 0 & 0 \\
0 & 0 & \frac{1}{\sqrt{2}} & 0 & \frac{1}{\sqrt{2}} &0 \\
0 & \frac{1}{\sqrt{2}} & 0 & -\frac{1}{\sqrt{2}} & 0 & 0 \\
0 & 0 & -\frac{1}{\sqrt{2}} & 0 &  \frac{1}{\sqrt{2}} & 0 \\
0 & \frac{1}{\sqrt{2}} & 0 & \frac{1}{\sqrt{2}} & 0 & 0 \\
0 & 0 & 0 & 0 & 0 & 1
\end{pmatrix}.
\]
\begin{ndefn}
	Given an $(n,n)$-clan $\gamma= c_1\cdots c_{2n}$, one obtains the \emph{base clan} associated to $\gamma$ by replacing each symbol $\tl{c}_i$ of the default signed clan $\tl\gamma$ by its signature.
\end{ndefn}
For example, the base clan for ${-}12334412{+}$ is ${-}{-}{-}{-}{+}{-}{+}{+}{+}{+}$. Now we can use the flags produced by Theorem \ref{thm:Dflag} to form the sects.

\begin{nrk}
A clan with no natural numbers is said to be \emph{matchless}. The base clan of a $DIII$ clan is clearly a matchless $DIII$ clan, and all matchless clans arise in this manner. Matchless clans correspond to closed orbits, which are also of minimum dimension. 
\end{nrk}

\begin{nprop} \label{lem:Dclanorbit} Let $Q_\gamma$ and $Q_\tau$ be $L$-orbits in $G/B$ corresponding to $DIII$ $(n,n)$-clans $\gamma$ and $\tau$. Then $Q_\gamma$ and $Q_\tau$ lie in the same $P$-orbit of $G/B$ if and only if $\gamma$ and $\tau$ have the same base clan.
\end{nprop}

\begin{proof} Assume that $\gamma$ has base clan $\tau$, where $\gamma=c_1\cdots c_{2n}$ and $\tau=t_1\cdots t_{2n}$, and let $V_\bullet(\gamma)=\la \mb{v}_1,\dots, \mb{v}_{2n} \ra$ and $V_\bullet(\tau)=\la \mb{u}_1, \dots, \mb{u}_{2n} \ra$ be the corresponding flags constructed by Theorem \ref{thm:Dflag}. As each clan has the same signature at symbols of the same index, they have the same default permutation. Then, we have two kinds of cases to examine.
	
	Suppose we have a family, $c_i,c_j,c_{2n+1-j},c_{2n+1-i}$. Then, Theorem \ref{thm:Dflag} will yield flag $V_\bullet(\gamma)$ with 
	\[(\mb{v}_i, \mb{v}_{j}) = (\frac{1}{\sqrt{2}}(\mb{e}_r+\mb{e}_{2n+1-s}),\frac{1}{\sqrt{2}}( \mb{e}_r -\mb{e}_{2n+1-s})) \]
	and
	\[(\mb{v}_{2n+1-j}, \mb{v}_{2n+1-i}) = (\frac{1}{\sqrt{2}}(-\mb{e}_s+\mb{e}_{2n+1-r}),\frac{1}{\sqrt{2}}( \mb{e}_s +\mb{e}_{2n+1-r})),\]
	where $n<r={\sigma}(i)$ and $n<s={\sigma}(2n+1-j)$. We also obtain the flag $V_\bullet(\tau)$ with 
	\[(\mb{u}_i, \mb{u}_{j})= (\mb{e}_r, \mb{e}_{2n+1-s}) \]
	and
	\[(\mb{u}_{2n+1-j}, \mb{u}_{2n+1-i}) = (\mb{e}_s,\mb{e}_{2n+1-r}).\]
	
	Then define a linear map by
	\begin{align*} 
	p^{r,s}:\ & \mb{e}_r \lmt \mb{e}_r+\mb{e}_{2n+1-s},  \\
	&\mb{e}_s \lmt \mb{e}_s-\mb{e}_{2n+1-r},  \\
	&\mb{e}_i\lmt \mb{e}_i \qquad\qquad\qquad\text{ for } i\neq r,s.
	\end{align*}
	It is again routine to check that this map defines an element of $P$, so $p^{r,s} \cdot V_\bullet(\tau)$ is a flag in the same $P$-orbit. Also, this map takes $\mb{u}_i$ to $\mb{v}_i$, and $\mb{u}_{2n+1-j}$ to the span of $\mb{v}_{2n+1-j}$, yielding pairs with the same span
	\[ (p^{r,s}\cdot \mb{u}_i,p^{r,s}\cdot \mb{u}_{j})\quad \text{and} \quad (\mb{v}_i, \mb{v}_{j}),\]
	and
	\[ (p^{r,s}\cdot \mb{u}_{2n+1-j},p^{r,s} \cdot \mb{u}_{2n+1-i})\quad \text{and} \quad (\mb{v}_{2n+1-j}, \mb{v}_{2n+1-i}).\]
	
	Now, after we act on the flag $V_\bullet(\tau)$ by the appropriate element of the form $p^{r,s}$ for each family,
	\[\{c_i=c_j, c_{2n+1-j}=c_{2n+1-i} \mid j\neq 2n+1-i\},\]
	then we obtain a flag which is an equivalent presentation of $V_\bullet(\gamma)$. Thus $Q_\gamma$ is in the same $P$-orbit as $Q_\tau$. 
	
	The proof of the converse is exactly as in type $CI$ case, which can be found in~\cite{aoSects}.
\end{proof}

By flipping the $L\bs G / B$ double cosets around and applying the map $\pi$, we obtain the following corollary. See also~\cite[Proposition 5.6]{aoSects}.

\begin{ncor} \label{projbase}
	Let $R_\gamma$ and $R_\tau$ be $B$-orbits in $G/L$ corresponding to clans $\gamma$ and $\tau$, and let \hbox{$\pi: G/L \to G/P$} denote the canonical projection. Then $\pi(R_\gamma)= \pi(R_\tau)$ if and only if $\gamma$ and $\tau$ have the same base clan.
\end{ncor}

Schubert cells of $SO_{2n}/P$ are in bijection with subsets $I\subs [2n]$ such that $\ab{I}=n$ and if $i \in I$, then $2n+1-i\not\in I$ \cite[p. 34]{billeyLak}. $P$ stabilizes the maximal isotropic subspace $E_n$, and in fact each $B$-orbit of $G/P$, denoted $C_I$, is represented by the isotropic subspace which is spanned by $\{\mb{e}_i \mid i \in I\}$. The subset $I$ can be associated to a matchless clan $\gamma_I$ by assigning
\begin{equation}\label{baseclanI} c_i =\begin{cases}
+ \quad \text{ if } i \in I \\
- \quad \text{ if } i \notin I,
\end{cases}\end{equation}
and just as in~\cite{aoSects}, $Bg_{\gamma_I}$P$= C_I$. Then we have the following analog of~\cite[Theorem 5.7]{aoSects}, whose proof is identical to the one given there, except for the fact that in this case $g_{\gamma_I}^{-1}=g_{\gamma_I}$, since it is the matrix of an even involution.

\begin{nthm} \label{decomp} Let $C_I$ be the Schubert cell corresponding to $I\subs [2n]$, and $\pi: G/L \to G/P$ the natural projection. Associate to $I$ a matchless clan $\gamma_I$ as in equation (\ref{baseclanI}), and denote the set of clans with base clan $\gamma_I$ by $\Sigma_I$. If $R_\gamma$ denotes the $B$-orbit of $G/L$ associated to the clan $\gamma$, then
	\begin{equation}
	\pi^{-1}(C_I) = \bigsqcup_{\gamma \in \Sigma_I} R_\gamma.
	\end{equation}
\end{nthm}
Consequently, each sect has a base clan which corresponds to a closed orbit, and the sects are in correspondence with Schubert cells. Further, each sect contains a unique maximal orbit, and the classes of closures of these orbits form a $\Z$-linear basis for the cohomology ring of $G/L$. We remark again that since $\pi: G/L\to G/P$ is an affine bundle with fibers isomorphic to $R_u(P)$, each sect can also be viewed as a decomposition of an affine space isomorphic to $C_I \x R_u(P)$ into $B$-orbits. %In the case of the unique dense Schubert cell, we study this decomposition in greater detail in Section~\ref{S:bigsect}.

\section{The weak order and its rank polynomial}\label{S:PartialOrder}
\subsection{The weak order on clans}

We continue with $B\subs SO_{2n}$ as the Borel subgroup of upper triangular matrices, and $L\cong GL_n$ as in (\ref{eq:PL}). Here we will describe the weak order on $DIII$ $(n,n)$-clans and calculate a recurrence for the rank polynomial of the weak order poset; see \cite{richardsonSpringer1990}, \cite{richardsonSpringer1993}, and \cite{wyserThesis} for further background and details on the weak order and its properties. We will first describe the geometric content of the weak order in terms of the $B$-orbits of $G/L$ (denoted $R_\gamma$ for clan $\gamma$) though it is of course equivalent to discuss $L$-orbits of $G/B$ or $B \x L$-orbits of $G$. 

Let $T\subs SO_{2n}$ be the maximal torus of diagonal matrices with Lie algebra $\mf{t}$. Note that this torus is $\theta$-stable and moreover is contained in $L=G^\theta$. By the condition defining the special orthogonal group, we have that
\[T=\{ \diag(t_1,\dots t_n, t_n^{-1},\dots,t_1^{-1})  \mid t_i \in \C^*) \} ,\]
so that 
\[\mf{t}=\{ \diag(a_1,\dots, a_n, -a_n, \dots,-a_1) \mid a_i \in \C\} .\]
We declare simple roots $\alpha_i := Y_i - Y_{i+1}$ for $1 \leq i \leq n-1$ and $\alpha_n:= Y_{n-1} + Y_n$ where $Y_i \in \mf{t}^* $ 
    is given by $Y_i(\diag(a_1,\dots, a_n,-a_n, \dots,-a_1))= a_i$. 

Corresponding to each simple root $\alpha_i$, there is a simple reflection $s_i$ in the Weyl group $W=N_G(T)/T$, and a minimal standard parabolic subgroup $P_{s_i}$. For any $DIII$ clan $\gamma$, $P_{s_i}\cdot R_\gamma$ contains a unique dense $B$-orbit $R_{\gamma'}$. To capture this relationship between $B$-orbits, we write $s_i \cdot R_\gamma = R_{\gamma'}$ and view this as an action of the set of simple reflections on the orbits. Note that if $\gamma\neq\gamma'$, then $R_{\gamma'}$ always has dimension equal to $\dim(R_\gamma)+1$.

Under this action, it is clear that $s_i\cdot(s_i \cdot R_\gamma)= s_i \cdot R_\gamma$ for any $i$ and $\gamma$. It is also true that the simple reflections obey the same braid relations when acting on $B$-orbits as they do in $W$ in its presentation as a Coxeter system. Thus, we actually have a monoid $M(W)$ which acts on the set of orbits and is generated by the simple reflections $\{s_i\}_{i=1}^n$ with relations $s_i^2=s_i$ plus braid relations. This monoid arises naturally in a degeneration of the Hecke algebra associated to $W$ as well \cite[Section 7]{richardsonSpringer1993}.

The \emph{weak order on DIII clans} is then defined as the transitive closure of the covering relations $\gamma \prec \gamma'$ whenever there is an $s_i$ such that $s_i\cdot R_\gamma=R_{\gamma'}$. In the following discussion, we may also abuse notation and write $s_i\cdot \gamma =\gamma'$ to mean the same. As weak order covering relations are labelled by simple reflections, the maximal chains of intervals in the weak order can be viewed as reduced expressions for elements of the orbit set, or alternatively for the underlying involutions in $W$, or their corresponding elements in $M(W)$. For more on this perspective (in type $A$ symmetric spaces), see \cite{canJoyceWyser}, \cite{hamMarbPawl1} and \cite{hamMarbPawl2}. 

Indeed, the simple reflections effectively act upon a clan $\gamma=c_1\cdots c_{2n}$ via its underlying involution $\sigma_\gamma$ through the following \emph{twisted action}. Consider $W$ as a subgroup of $\mc{S}_{2n}$, so that we have $s_i=(i\ (i+1))((2n-i)(2n+1-i))$ for $1\leq i \leq n-1$ and $s_n= (n\ (n+2))((n-1)(n+1))$. 

\begin{nprop} \label{prop:twist}
Suppose $s_i\sigma_\gamma$ is longer than $\sigma_\gamma$ as an element of $W$ for the $DIII$ clan $\gamma=c_1\cdots c_{2n}$. Then 
\begin{enumerate} 
\item if $s_i \sigma_\gamma s_i \neq \sigma_\gamma$, then $ s_i \cdot \gamma$ is the permutation action of $s_i$ on the symbols of $\gamma$ which results in underlying involution $\sigma_{\gamma'}=s_i \sigma_\gamma s_i$;
\item if $s_i \sigma_\gamma s_i = \sigma_\gamma$ for any $1\leq i \leq n-1$, and $c_i$ and $c_{i+1}$ are opposite signs $+/-$ (in either order), then $s_i \cdot \gamma$ replaces the appropriate signed fixed points by natural number pairs to achieve modified underlying involution $\sigma_{\gamma'}=s_i \sigma_\gamma$;
\item if $c_{n-1}c_nc_{n+1}c_{n+2}={+}{+}{-}{-}$ or ${-}{-}{+}{+}$, then $s_i \cdot \gamma$ replaces these symbols by the pattern $1212$, resulting in underlying involution $\sigma_{\gamma'}=s_i\sigma_\gamma$.
\end{enumerate}
Otherwise $s_i\cdot \gamma =\gamma$.
\end{nprop}
\begin{proof}
This is evident from \cite[Theorem 5.4.1]{richardsonSpringer1993}, which also applies to the $AIII$ and $CI$ cases. Full translation of the action of simple reflections into $DIII$ clan notation is also worked out (with examples) in \cite[Section 5.2.2]{wyserThesis}.
\end{proof}

\begin{nexap}
    For the $(4,4)$-clan $\gamma={+}{-}1122{+}{-}$, from the rules above we obtain $s_1\cdot \gamma = 11223344$ and $s_2\cdot \gamma = {+} 1 {-} 12{+} 2{-}$, while $s_3\cdot \gamma = \gamma$ and $s_4\cdot \gamma =\gamma$. See also Figure~\ref{fig:del4}.
\end{nexap}
From this discussion, it follows that the weak order poset on $DIII$ $(n,n)$-clans, denoted $(\del(n),\prec)$, is ranked (graded) by the length of underlying involutions in terms of the twisted action indicated by Proposition~\ref{prop:twist}. Note that this is often different than the the length of the underlying involution as a Weyl group element.

\begin{ndefn}
We define the length $L(\gamma)$ of a $DIII$ $(n,n)$-clan $\gamma$ as the length of its underlying involution $\sigma_\gamma$ under the twisted action.
\end{ndefn}
See \cite[Section 5]{richardsonSpringer1990} for various properties of the twisted action and this length function. As an example, matchless clans all have the identity as underlying involution, and so they have length $0$. 

\begin{landscape}
	\begin{figure}[h]
		\begin{center}
			\begin{tikzpicture}[scale=.22]
			\node at (-47,0) (a1) {$\blue{{+}{+}{+}{+}}$};
			\node at (-34,0) (a2) {$\blue{{+}{+}{-}{-}}$};
			\node at (-21,0) (a3) {$\blue{{+}{-}{+}{-}}$};
			\node at (-8,0) (a4) {$\blue{{+}{-}{-}{+}}$};
			\node at (5,0) (a5) {$\blue{{-}{+}{+}{-}}$};
			\node at (18,0) (a6) {$\blue{{-}{+}{-}{+}}$};
			\node at (33,0) (a7) {$\blue{{-}{-}{+}{+}}$};
			\node at (47,0) (a8) {$\blue{{-}{-}{-}{-}}$};

			\node at (-47,10) (b1) {${+}{+}1212{-}{-}$};
			\node at (-34,10) (b2) {${+}{-}1122{+}{-} $};
			\node at (-21,10) (b3) {${+}11{-}{+}22{-}$};
			\node at (-8,10) (b4) {$11{+}{-}{+}{-}22$};
			\node at (5,10) (b5) {$11{-}{+}{-}{+}22$};
			\node at (18,10) (b6) {${-}11{+}{-}22{+}$};
			\node at (33,10) (b7) {${-}{+}1122{-}{+}$};
			\node at (47,10) (b8) {${-}{-}1212{+}{+}$};

			\node at (-41,20) (c1) {${+}1{+}21{-}2{-}$};
			\node at (-26,20) (c2) {${+}1{-}12{+}2{-}$};
			\node at (-11,20) (c3) {$1{+}1{-}{+}2{-}2$};
			\node at (0,20) (c4) {$11223344$};
			\node at (11,20) (c5) {$1{-}1{+}{-}2{+}2$};
			\node at (26,20) (c6) {${-}1{+}12{-}2{+}$};
			\node at (41,20) (c7) {${-}1{-}21{+}2{+}$};

			\node at (-41,30) (d1) {${+}12{+}{-}12{-}$};
			\node at (-26,30) (d2) {$1{+}{+}21{-}{-}2	$};
			\node at (-11,30) (d3) {$1{+}{-}12{+}{-}2$};
			\node at (0,30) (d4) {$12123434$};
			\node at (11,30) (d5) {$1{-}{+}12{-}{+}2$};
			\node at (26,30) (d6) {$1{-}{-}21{+}{+}2$};
			\node at (41,30) (d7) {${-}12{-}{+}12{+}$};

			\node at (-35,40) (e1) {$1{+}2{+}{-}1{-}2$};
			\node at (-11,40) (e2) {$12213443$};
			\node at (11,40) (e3) {$12341234$};
			\node at (35,40) (e4) {$1{-}2{-}{+}1{+}2$};

			\node at (-22,50) (f1) {$12{+}{+}{-}{-}12$};
			\node at (0,50) (f2) {$12342143$};
			\node at (22,50) (f3) {$12{-}{-}{+}{+}12$};
			
			\node at (0,60) (g1) {$12343412$};

			\path 
			(a1) edge [thick] node[left]  {\tiny$4$} (b1)
	        (a2) edge [thick] node[left,pos=.2]  {\tiny$4$} (b1)
			(a2) edge [thick] node[left,pos=.2]  {\tiny$2$} (b3)
			(a3) edge [thick] node[right,pos=.2]  {\tiny$3$} (b2)
			(a3) edge [thick] node[left,pos=.2]  {\tiny$2$} (b3)
			(a3) edge [thick] node[left,pos=.2]  {\tiny$1$} (b4)
		    (a4) edge [thick] node[right, pos=.2]  {\tiny$3$} (b2)
			(a4) edge [thick] node[left,pos=.2]  {\tiny$1$} (b5)
		    (a5) edge [thick] node[right,pos=.2]  {\tiny$1$} (b4)
			(a5) edge [thick] node[left,pos=.2]  {\tiny$3$} (b7)
			(a6) edge [thick] node[right,pos=.2]  {\tiny$1$} (b5)
			(a6) edge [thick] node[left,pos=.2]  {\tiny$2$} (b6)
		    (a6) edge [thick] node[left,pos=.2]  {\tiny$3$} (b7)
			(a7) edge [thick] node[left,pos=.2]  {\tiny$2$} (b6)
		    (a7) edge [thick] node[left]  {\tiny$4$} (b8)
			(a8) edge [thick] node[left]  {\tiny$4$} (b8)
			
			 (b1) edge [thick] node[left]  {\tiny$2$} (c1)
			(b2) edge [thick] node[left,pos=.2]  {\tiny$2$} (c2)
			(b2) edge [thick] node[left,pos=.5]  {\tiny$1$} (c4)
			(b3) edge [thick] node[right,pos=.7]  {\tiny$1$} (c1)
			(b3) edge [thick] node[left]  {\tiny$3$} (c2)
		 (b3) edge [thick] node[right,pos=.3]  {\tiny$1$} (c3)
			(b4) edge [thick] node[left]  {\tiny$2$} (c3)
			(b4) edge [thick] node[left]  {\tiny$3$} (c4)
			 (b5) edge [thick] node[left]  {\tiny$3$} (c4)
			 (b5) edge [thick] node[left]  {\tiny$2$} (c5)
			(b6) edge [thick] node[left,pos=.2]  {\tiny$3$} (c5)
			 (b6) edge [thick] node[right,pos=.7]  {\tiny$3$} (c6)
			(b6) edge [thick] node[left,pos=.7]  {\tiny$3$} (c7)
			(b7) edge [thick] node[right,pos=.5]  {\tiny$1$} (c4)
			 (b7) edge [thick] node[left,pos=.7]  {\tiny$2$} (c6)
			(b8) edge [thick] node[left]  {\tiny$2$} (c7)

			(c1) edge [thick] node[left]  {\tiny$3$} (d1)
			(c1) edge [thick] node[left,pos=.2]  {\tiny$1$} (d2)
			(c2) edge [thick] node[right,pos=.2]  {\tiny$4$} (d1)	
			(c2) edge [thick] node[left,pos=.2]  {\tiny$1$} (d3)
			(c3) edge [thick] node[right,pos=.2]  {\tiny$4$} (d2)
		    (c3) edge [thick] node[left]  {\tiny$3$} (d3)
			(c4) edge [thick] node[left]  {\tiny$2$} (d4)
			%(c4) edge [thick] node[left]  {\tiny$3$} (d4)
			(c5) edge [thick] node[left]  {\tiny$3$} (d5)
			(c5) edge [thick] node[left,pos=.2]  {\tiny$4$} (d6)
			(c6) edge [thick] node[right,pos=.2]  {\tiny$1$} (d5)
			(c6) edge [thick] node[left,pos=.2]  {\tiny$4$} (d7)
		    (c7) edge [thick] node[right,pos=.2]  {\tiny$1$} (d6)
		    (c7) edge [thick] node[left]  {\tiny$3$} (d7)
			
			(d1) edge [thick] node[left]  {\tiny$1$} (e1)
		    (d2) edge [thick] node[left]  {\tiny$3$} (e1)
			(d3) edge [thick] node[left]  {\tiny$2$} (e2)
			(d3) edge [thick] node[left]  {\tiny$4$} (e1)
			(d4) edge [thick] node[left]  {\tiny$1$} (e2)
			(d4) edge [thick] node[left,pos=.2]  {\tiny$4$} (e3)
			(d5) edge [thick] node[right,pos=.2]  {\tiny$2$} (e2)
			(d5) edge [thick] node[left]  {\tiny$4$} (e4)
			(d6) edge [thick] node[left]  {\tiny$3$} (e4)
			(d7) edge [thick] node[left]  {\tiny$1$} (e4)

			(e1) edge [thick] node[left]  {\tiny$2$} (f1)
			(e2) edge [thick] node[left]  {\tiny$4$} (f2)
			(e3) edge [thick] node[left]  {\tiny$3$} (f2)
			(e4) edge [thick] node[left]  {\tiny$2$} (f3)

			(f1) edge [thick] node[left]  {\tiny$4$} (g1)
			(f2) edge [thick] node[left]  {\tiny$2$} (g1)
			(f3) edge [thick] node[left]  {\tiny$4$} (g1)
			
			(a2) edge [thick, red]    node[right] {}  (b2)
			(a4) edge [thick, red]    node[right] {}  (b4)
			(a5) edge [thick, red]    node[right] {}  (b5)
			(a7) edge [thick, red]    node[right] {}  (b7)
			(c3) edge [thick, red]    node[right] {}  (d4)
			(c4) edge [thick, red]    node[right] {}  (d3)
			(c4) edge [thick, red]    node[right] {}  (d5)
			(c5) edge [thick, red]    node[right] {}  (d4)
			(d1) edge [thick, red]    node[right] {}  (e3)
			(d2) edge [thick, red]    node[right] {}  (e2)
			(d2) edge [thick, red]    node[right] {}  (e3)
			(d6) edge [thick, red]    node[right] {}  (e2)
			(d6) edge [thick, red]    node[right] {}  (e3)
			(d7) edge [thick, red]    node[right] {}  (e3)
			(e2) edge [thick, red]    node[right] {}  (f1)
			(e2) edge [thick, red]    node[right] {}  (f3)
			(e4) edge [thick, red]    node[right] {}  (f2);

			\end{tikzpicture}
			
			\caption{The weak and full closure orders on $\del(4)$. }
			\label{fig:del4}
		\end{center}
	\end{figure}
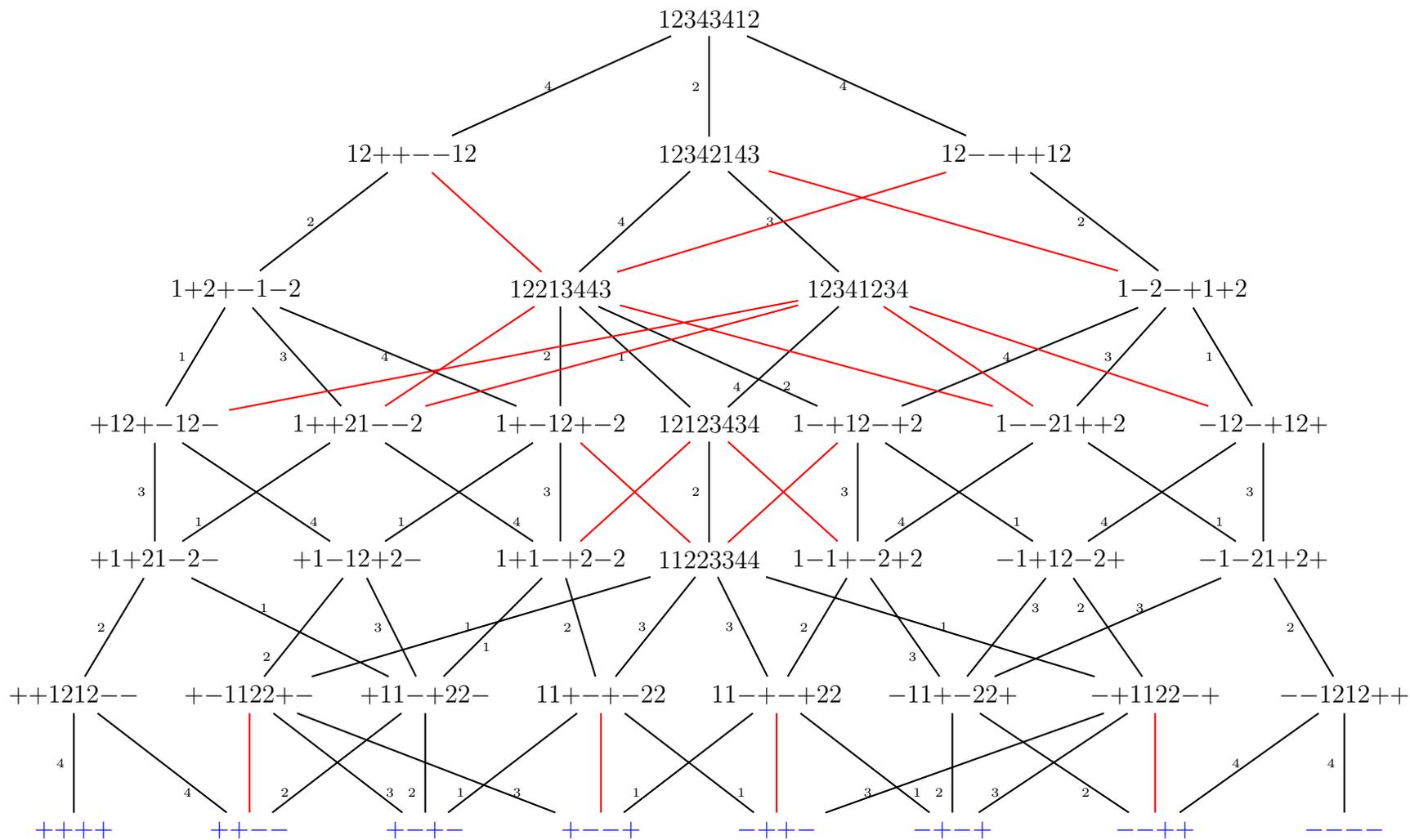
\end{landscape}

\subsection{Rank polynomial} \label{ssec:lgf}

We will write $A_n(t)$ to denote the \emph{rank polynomial} of the weak order poset. That is, $A_n(t)$ is the polynomial in $t$ for which the coefficient of $t^k$ is equal to the number of $DIII$ clans of length $k$. We may also call $A_n(t)$ the \emph{length generating function} for $DIII$ $(n,n)$-clans. In order to compute a recurrence for $A_n(t)$, we shall make use of a formula for the length of a clan given purely in terms of the string $\gamma=c_1\cdots c_{2n}$. 

First, we need some auxiliary notation. We will partition the natural number pairs of $\gamma$ into two sets. Let 
\[ \Pi_0:= \{ (c_i,c_j) \mid c_i=c_j\in \N \text{ and } 1\leq i \leq n < j \leq 2n \},\]
and
\[ \Pi_1:= \{ (c_i,c_j) \mid c_i=c_j\in \N \text{ and } 1\leq i < j \leq n \text{ or } n+1\leq i< j \leq 2n\}.\]
If a pair of mates $(c_i, c_j)$ in either one of these sets, then its opposing pair is in the same set. Then we can write $\ab{\Pi_0}= 2z$ and $\ab{\Pi_1}= 2y$ for integers $z$ and $y$.

For a natural number $a=c_i=c_j$ which appears in $\gamma=c_1\cdots c_{2n}$, the \emph{spread} of $a$ is defined as the quantity \mbox{$s(a):=j-i$}. The \emph{weave} of $a$ will be the quantity 
\[w(a):=\#\{b\in \N \mid b = c_u=c_t \text{ and } u<i<t<j\}.\]

\begin{nprop}
In the notation above, the length of a clan $\gamma$ is equal to
\begin{equation}\label{eq:Length}
L(\gamma)=\frac{1}{2}\lt(\bigg(\sum_{a=c_i=c_j} s(a)-w(a)\bigg) -z\rt).
\end{equation}
\end{nprop}
\begin{proof}
    The reader may verify that the formula is equivalent to the one that appears at \cite[p. 2724]{mcgovernTrapa}, after subtracting the dimension of a closed orbit. %Further, this formula is equivalent to the $D$-inversion number of the underlying involution $\sigma_\gamma$, which also calculates the length of $\sigma_\gamma$ as in \cite[Section 8.2]{bjornerBrenti}.
\end{proof}
\begin{nexap}
    Take the $DIII$ clan $\gamma={+}{+}1212{-}{-}$. Each pair of mates contributes a spread of 2, but only $2=c_4=c_6$ has a weave of 1, and together these pairs are the only $\Pi_0$ family so $z=1$. Thus, $L(\gamma)=1$; see Figure~\ref{fig:del4}.
\end{nexap}
	\begin{comment}
	We define the length $L(\gamma)$ of a $DIII$ $(n,n)$-clan $\gamma= c_1 c_2 \dots c_{2n}$ as follows:
	\begin{equation*}
	L(\gamma) = \frac{1}{2} \biggr(l(\gamma) - \#\{a\in \N \;|\; a=c_u =c_t \; \text{and}\; u\leq n<t\leq 2n+1-u \} \biggr)
	\end{equation*}
	where $l(\gamma) = \sum\limits_{\substack{c_i=c_j \in \N\\ i<j}} \big(j-i-\# \{a\in \N \;|\; a=c_u =c_t \; \text{and}\; u< i<t< j\}\big).$
 Note that $l$ is the length function for $(p,q)$-clans (type $AIII$) which appears in~\cite{yamamoto}. Moreover, in the notation of Section~\ref{S:Formula}, 
\begin{equation}
L(\gamma)= \frac{1}{2}(l(\gamma)-p)
\end{equation}
where $2p=\ab{\Pi_0}$.  
\end{comment}

\begin{nrk}
The expression $\sum_{a=c_i=c_j} s(a)-w(a)$ computes the length for $(n,n)$-clans in type $AIII$, as appears in~\cite{yamamoto}.
\end{nrk}

Note that the inclusion poset of Borel orbit closures in $SO_{2n}/GL_n$ contains all of the weak order relations on $DIII$ $(n,n)$-clans, possibly plus some additional relations.\footnote{The order relation of this poset is often called the \emph{(full) closure order} or \emph{Bruhat order}.} Thus, the length function also provides a grading of the closure containment poset. It follows from the description of the weak order that the unique maximal clan in both posets is of the form
\begin{eqnarray} \label{eq:gamma0}
\gamma_0 &=& 1 2 \dots (n-1) n (n-1) n \dots 1 2 \quad \text{if}\; n \; \text{is even, and} \nonumber\\  \\
\gamma_0 &=& 1 2 \dots (n-1) n {+} {-} (n-1) n \dots 1 2 \quad \text{if}\; n \; \text{is odd}.\nonumber
\end{eqnarray}
\begin{nrk}
	If $\gamma $ is the $DIII$ clan corresponding to the Borel orbit $R_\gamma$, then the dimension of $R_\gamma$ is equal to $L(\gamma) + c$, where $c$ is the dimension of any closed Borel orbit in $SO_{2n}/GL_n$. Thus, studying $L(\gamma)$ is equivalent to studying dimensions of Borel orbits in the type $DIII$ symmetric space (or $GL_n$-orbits in $SO_{2n}/B$). Moreover, the dimension of all closed orbits is equal to the dimension of the flag variety of $GL_n$ (or of $R_u(P)$), which is $\frac{n(n-1)}{2}$.
\end{nrk}

\begin{nprop}
The length of $\gamma_0$, the maximal element in the weak order poset, is $\frac{n(n-1)}{2}$.
\end{nprop}
\begin{proof}
$SO_{2n}$ has the dimension of its Lie algebra, which consists of skew-symmetric $2n \x 2n$ matrices. This is $\frac{2n(2n-1)}{2}$ dimensional. $GL_n$ has dimension $n^2$, and so $$\dim SO_{2n}/GL_n= n(2n-1)-n^2=n(n-1).$$ The maximal element corresponds to a dense orbit $R_{\gamma_0}\subs SO_{2n}/GL_n$, so $\dim R_{\gamma_0}=n(n-1)$ as well. By the preceding remark, this is also equal to $L(\gamma_0)+\frac{n(n-1)}{2}$, finishing the proof.
\end{proof}

As before, $\del(n)$ denotes the set of $DIII$ $(n,n)$-clans. We now reintroduce $A_n(t)$.
%Note that if $\gamma$ is a $DIII$ clan then $\Flip(\gamma)$ is not, though we will still act on it combinatorially by the simple reflections as if it were.
\begin{ndefn}
The length generating function of $\Delta (n)$ is defined by $$A_n(t) = \sum\limits_{\gamma \in \Delta(n)} t^{L(\gamma)}. $$
\end{ndefn}
We also define the \emph{flip} of $\gamma=c_1\cdots c_nc_{n+1} \cdots c_{2n}$ by
$$
\text{Flip}(\gamma):=c_1\cdots c_{n+1}c_n \cdots c_{2n}.
$$
Now, we provide a recurrence relation for $A_n(t)$.
\begin{nthm} \label{thm:lgf}
The length generating function $A_n(t)$ satisfies the following recurrence for $n\geq 3$:
\begin{equation} \label{eq:lenrec}
A_n(t)=2A_{n-1}(t) + (t+ t^2 + \cdots t^{n-2}+2t^{n-1}+t^n +\cdots +t^{2n-3} )A_{n-2}(t).
\end{equation}
\end{nthm}
\begin{proof}
We break this into two parts, one for each of the recursive terms. 

\textbf{Coefficient of }$\mb{A_{n-1}(t)}:$ Let $\gamma$ be an arbitrary clan from $\del(n-1)$. Then, we can always create a new clan $+\gamma-\in\del(n)$ simply by inserting a $+$ at the beginning of the string and appending a $-$ at the end of the string. It is clear that this procedure does not affect the value of the length function. 

We can create a different clan $\Flip(-\gamma+)\in\del(n)$ similarly, where the flip is required to ensure that there are an even number of $-$'s and $\Pi_1$ pairs among the first $n$ symbols. In this situation, there are a few possibilities. Let $\gamma=c_2\cdots c_n c_{n+1} \cdots c_{2n-1}$ for convenience.

$\mb{[c_nc_{n+1}=\pm\mp]}:$ Attaching the new symbols and flipping has no consequence for any component of the length function computation of (\ref{eq:Length}).

$\mb{[\gamma = \cdots a \cdots ab\cdots b \cdots]}:$ Attaching $-$ and $+$ has no effect. Upon flipping, $s(a)$ and $s(b)$ both increase by one, but so does $w(b)$ and $z$, so there is no net effect on the length function.

$\mb{[\gamma = \cdots a \cdots ba\cdots b \cdots]}:$ Identical to the previous case, but change ``increase'' to ``decrease.''

We see that given an arbitrary $(n-1, n-1)$-clan, we can create two different $(n,n)$-clans for which the length function evaluates the same, accounting for the first term in the equation (\ref{eq:lenrec}). These comprise all of the clans in $\del(n)$ which start and end with $+$ or $-$.

\textbf{Coefficient of }$\mb{A_{n-2}(t)}:$ Now let $\gamma$ be an arbitrary clan from $\del(n-2)$. We obtain a new clan $\gamma'=c'_1\cdots c'_{2n}\in \del(n)$ by inserting $a\in \N$ as $c'_1$ and $c'_{i}$ and $b\in \N$ as $c'_{2n+1-i}$ and $c'_{2n}$. Observe that $a$ and $b$ each contribute a spread of $i-1$, and $w(a)=0$ for any choice $i$.

If $i\leq n$, then both $a$ and $b$ go in as $\Pi_1$ pairs, so $z$ is unchanged. If $w(b)$ results positive due to any natural number pair $(c'_u,c'_t)$, this contribution will cancel in the length formula by the fact that the first $b$ at $c'_{2n+1-i}$ increases the spread of that pair by one, compared to its placement in $\gamma$. The insertion of $b$ cannot affect the weave of any other natural number because the last symbol is $b$. If $a$ contributes to the weave of any other natural number $r$, this is likewise cancelled out by an increase of one in $s(r)$. Thus, the length is only affected by the spreads of $a$ and $b$, and increases by $i-1$ on balance. Since this holds for any choice of $2\leq i\leq n$, we see that $(t+t^2+ \dots t^{n-1}) A_{n-2}(t)$ appears in the recurrence. 

Now suppose $i>n$. Both $a$ and $b$ go in as $\Pi_0$ pairs, so $z$ increases by one. As with the previous case, weave contributions of $a$ and $b$ cancel with spread contributions to other numbers with one exception: the pair of $a$'s $(c'_1,c'_i)$ contributes one to the weave of $b$ which is not compensated for in any other manner. Thus, compared to the length of $\gamma$, $L(\gamma')$ is increased by $i-1$ from the spreads of $a$ and $b$, and but diminished by one from the change in $z$ and $w(b)$. In all, each choice of $n<i\leq 2n-1$ gives a different clan whose length is $i-2$ greater than $L(\gamma)$, accounting for an additional term of
$(t^{n-1} +t^{n}+\dots t^{2n-3})A_{n-2}(t)$ in the recurrence formula. 

Adding these cases together gives the claim.  
\end{proof}
For consideration, we mention that $A_1(t)=1$, $A_2(t)=t+2$, and, in view of Figure~\ref{fig:del4}, $A_4(t)=t^6+3t^5+4t^4+7t^3+7t^2+8t+8$. In that figure, the black edges between $\gamma\prec \gamma'$ are labelled with the index $i$  such that $s_{i} \cdot \gamma = \gamma'$, while the red edges are those from the full closure order on corresponding orbits (see \cite[Proposition 4.2]{richardsonSpringer1993} for how to obtain the the full closure order from the weak order). The closed orbits in blue are represented by just their first four symbols due to space considerations.

 The recurrence for $A_n(t)$ yields the following statements about $\del_n$, the number of $DIII$ $(n,n)$-clans. In the next section, we will give an explicit formula for $\del_n$.

\begin{ncor}
	For all $n\geq 3$, the number of $DIII$ clans satisfies the recurrence relation 
	\begin{equation} \label{eq:rec}
	\del_n =2\del_{n-1}+(2n-2)\del_{n-2},
	\end{equation}
	and assigning $\del_0=1$, $\del_n$ has exponential generating function 
	\begin{equation}\label{eq:dgen}
	\sum_{n=0}^\infty \del_n\frac{x^n}{n!} = \frac{1}{2}(e^{2x+x^2}+1).
	\end{equation}
\end{ncor}
\begin{proof} The recurrence follows by substituting $t=1$ into equation~(\ref{eq:lenrec}). One can check that $y=\frac{1}{2}e^{2x+x^2}$ solves the relevant second-order linear homogeneous ordinary differential equation, $y''-2(x+1)y'+2y=0$, and the addition of 1 is just to satisfy the initial condition $y(0)=1$ coming from $\del_0=1$.
\end{proof}

\section{Bijective combinatorics for $DIII$ clans}\label{S:Formula}
\subsection{A Formula}
Before exploring bijections of clans with other combinatorial families, we will record an explicit formula for the number of $DIII$ $(n,n)$-clans $\del_n$. %As before, $\del(n)$ denotes the set of $DIII$ $(n,n)$-clans of and $\del_n$ its cardinality.  
%\begin{nlem} If $\gamma=c_1 \cdots c_{2n}$ is an $(n,n)$-clan with $k$ pairs of natural numbers among its symbols, then $k$ is even.
%\end{nlem}
%\begin{proof} Suppose $c_i=c_j =a\in \N$. By property 2 of type $DIII$ clans, $j\neq 2n+1-i$. Since $\gamma$ is skew-symmetric, upon taking the negative transpose we have that $a\neq c_{2n+1-i}=c_{2n+1-j}=b \in \N$. Thus, each natural number pair in $\gamma$ has an ``opposing'' natural number pair with which it exchanges place upon reversal. Since natural number pairs then come in twos (so accounting for four symbols $c_i$ in total), $k$ is even.
%\end{proof}

Recall from Section~\ref{ssec:lgf} that for any $DIII$ clan, the sets $\Pi_0$ and $\Pi_1$ both have even cardinality, as each pair of mates is half of a family where opposing pairs lie in the same set. Let $\delta_{r,n}$ denote the number of $DIII$ clans which contain $r$ families, or equivalently $2r$ pairs of mates. %Next, observe that if $\gamma$ contains $k$ pairs of natural numbers, then exactly $k$ of the symbols among $c_1\cdots c_n$ will be natural numbers, by the skew-symmetry condition.In order to count $\Delta_n$, first we will count $\delta_{r,n}$. 
Throughout this section, we make use of the fact that a $DIII$ $(n,n)$-clan is determined by the symbols of its ``first half'' $c_1\cdots c_n$, plus the knowledge of which $\Pi_0$ pairs are opposing.

\begin{nlem} With the above notation,
	\begin{equation} \label{eq:delknLem}
	\delta_{r,n} = 2^{n-2r-1} \binom{n}{2r} \frac{2r!}{r!}.
	\end{equation}
\end{nlem}
\begin{proof} There are $\binom{n}{2r}$ choices for where to place natural numbers among $c_1\cdots c_n$. There are $(2r-1)!!=\frac{2r!}{r!2^r}$ ways to form $r$ pairs from these, and $2^r$ ways to decide which of these pairs are in $\Pi_1$ and which are first mates of distinct opposing $\Pi_0$ pairs. Then there are $2^{n-2r-1}$ ways to place $\pm$ symbols at the remaining entries with appropriate parity to satisfy condition 3 of Definition~\ref{defn:Dclans}. Multiply.
\end{proof}
The following results immediately from summing over possible values of $r$.
\begin{nprop}\label{cor:Dformula}
	The number of $DIII$ $(n,n)$-clans  is 
	\begin{equation}
	\label{eq:Deln}
	\del_n= \sum_{r=0}^{\floor{\frac{n}{2}}} 2^{n-2r-1}\frac{n!}{r!(n-2r)!}.
	\end{equation}
\end{nprop}
%\begin{proof} The formula comes from rewriting the equation (\ref{eq:delknLem}) in terms of $r$, expanding the binomial coefficient, and then summing over possible values for $r$. 
%\end{proof}

The first values of $\del_n$, beginning with $n=1$, are 1, 3, 10, 38, 156, 692, 3256,\dots   This is in fact the number of inequivalent placements of $2n$ rooks on a $2n \x 2n$ board having symmetry across each diagonal \cite[\href{http://oeis.org/A000902}{A000902}]{OEIS}, as we will show next.

\subsection{Rooks and Pyramids} \label{subsec:rookpyr}
\begin{figure}[!htb] 
	\begin{center}
		\begin{tikzpicture}[scale=0.5]
		\draw  (0,0) grid +(6,6);
		\draw[thick] 
		(0,0) -- (6,0) -- (6,6) -- (0,6) -- (0,0);
		\draw[thin, color=red!50] (0,0)--(6,6) 
		(6,0) -- (0,6);
		\draw (0,6) node[anchor=east] {$d$} 
		(6,6) node[anchor=west] {$d'$};
		\filldraw[brown!50!black] (0.5, 1.5) circle (8pt)
		(1.5, 0.5) circle (8pt)
		(2.5, 2.5) circle (8pt)
		(3.5, 3.5) circle (8pt)
		(4.5, 5.5) circle (8pt)
		(5.5, 4.5) circle (8pt);
		\draw[thin, color=red!50] (0,0)--(6,6) 
		(6,0) -- (0,6);
		\draw (0,6) node[anchor=east] {$d$} 
		(6,6) node[anchor=west] {$d'$};
		\end{tikzpicture}
		\hspace{2cm}
		\begin{tikzpicture}[scale=0.5]
		\draw  (0,0) grid +(6,6);
		\draw[thick] 
		(0,0) -- (6,0) -- (6,6) -- (0,6) -- (0,0);
		\filldraw[brown!50!black] (0.5, 4.5) circle (8pt)
		(1.5, 5.5) circle (8pt)
		(2.5, 3.5) circle (8pt)
		(3.5, 2.5) circle (8pt)
		(4.5, 0.5) circle (8pt)
		(5.5, 1.5) circle (8pt);
		\end{tikzpicture}
		\caption{A diagonally invariant $6 \x 6$ rook placement and its rotational equivalent.}
		\label{boards1}
	\end{center}
\end{figure}
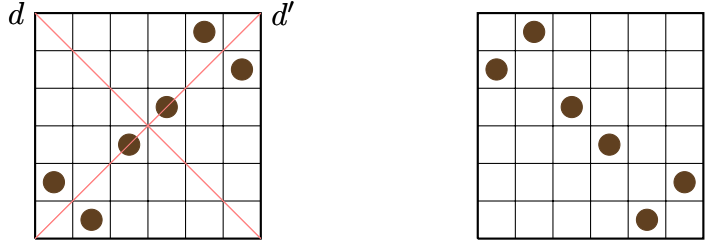

The rook problem asks how many ways $n$ rooks can be placed on on an $n \x n$ board so that none can attack any other. Necessarily, each placement will exhibit exactly one rook in each row and each column so rook placements correspond to permutation matrices in an obvious way, giving the solution of $n!$. In~\cite{lucas}, Lucas refines this question to ask how many placements possess symmetry with respect to a given subgroup of the dihedral group $D_8$, which acts as the symmetry group of the board. We are interested in rook placements which are invariant under reflection across each main diagonal $d$ and $d'$ as depicted in Figure~\ref{boards1}.

Furthermore, we are only interested in these placements up to equivalence, where two placements are said to be \emph{equivalent} if there is an element of $D_8$ that transforms one to the other. The following statements on rook placements are given without proof in~\cite{lucas}, but we provide brief proofs for completeness.

\begin{nprop}\label{prop:d8}
	When $n\geq 2$, there is no placement of $n$ rooks on an $n \x n$ board with symmetry group $D_8$.
\end{nprop}

\begin{proof}
       Recall that the dihedral group $D_8$ is generated by reflection $F$ across the diagonal $d$ and counter clockwise rotation by $\frac{\pi}{2}$ denoted $R$, and we have $R^4 = F^2 = (RF)^2 = e$ where $e$ is the identity element. It suffices to show that a rook placement cannot be symmetric with respect to both $F$ and $R$ unless $n=1$.

Consider a rook placement as an $n \x n$ permutation matrix $v$. Let $w_0$ be the permutation matrix with 1's along the antidiagonal $d'$ and 0's elsewhere. Then $R\cdot v= w_0v^{-1}$ and since the transpose of a permutation matrix is its inverse, $F\cdot v = v^{-1}$. Hence, if $v$ is invariant under both $R$ and $F$, we have 
 \[v=w_0v^{-1} \quad \text{and} \quad v=v^{-1},\]
 which implies $v^2=w_0$ and $v^2=I_n$, the $n \x  n$ identity matrix. This is clearly impossible unless $n=1$.
\end{proof}
\begin{ncor} When $n\geq 2$, every placement of $n$ rooks on an $n\x n$ board that is symmetric with respect to reflection across both diagonals is equivalent to exactly one other arrangement.
\end{ncor}
\begin{proof} Note that the reflection across $d'$ is given by $R^{-1}FR$. Together, $F$ and $R^{-1}FR$ generate a subgroup $V$ of $D_8$ isomorphic to the Klein four-group $\Z_2 \x \Z_2$. $V$ has index $2$ in $D_8$, so it is a maximal proper subgroup. Then, by the previous proposition, any rook placement stabilized by $V$ has exactly $V$ as its symmetry group, so by the orbit-stabilizer theorem, the size of the orbit of a diagonally symmetric rook placement under the action of $D_8$ is just
	\[ [D_8 : V] = \frac{\ab{D_8}}{\ab{V}}= \frac{8}{4}=2.\] 
\end{proof}
\begin{nrk}
For a diagonally symmetric rook placement $v$, $R\cdot v$ is necessarily different than (though equivalent to) $v$, so it is the other element of the orbit in the preceding corollary.
\end{nrk}
\begin{figure}[!t] 
	\begin{center}
		\begin{tikzpicture}[scale=0.5]
		\draw  (0,0) grid +(4,4);
		\draw[thick] 
		(0,0) -- (4,0) -- (4,4) -- (0,4) -- (0,0);
		\filldraw[brown!50!black]
		(0.5, 3.5) circle (8pt)
		(1.5, 2.5) circle (8pt)
		(2.5, 1.5) circle (8pt)
		(3.5, 0.5) circle (8pt);
		\draw[thin, dashed, color=red!50] (0,0)--(4,4) 
		(4,0) -- (0,4);
		\draw[<->] (2,-0.5) to (2,-1.5);					    
		\begin{scope}[xshift=8cm]
		\draw  (0,0) grid +(4,4);
		\draw[thick] 
		(0,0) -- (4,0) -- (4,4) -- (0,4) -- (0,0);
		\filldraw[brown!50!black]
		(0.5, 2.5) circle (8pt)
		(1.5, 3.5) circle (8pt)
		(2.5, 0.5) circle (8pt)
		(3.5, 1.5) circle (8pt);
		\draw[thin, dashed, color=red!50] (0,0)--(4,4) 
		(4,0) -- (0,4);
		\draw[<->] (2,-0.5) to (2,-1.5);	
		\end{scope}
		\begin{scope}[xshift=16cm]
		\draw  (0,0) grid +(4,4);
		\draw[thick] 
		(0,0) -- (4,0) -- (4,4) -- (0,4) -- (0,0);
		\filldraw[brown!50!black]
		(0.5, 3.5) circle (8pt)
		(1.5, 1.5) circle (8pt)
		(2.5, 2.5) circle (8pt)
		(3.5, 0.5) circle (8pt);
		\draw[thin, dashed, color=red!50] (0,0)--(4,4) 
		(4,0) -- (0,4);
		\draw[<->] (2,-0.5) to (2,-1.5);						    
		\end{scope}
		
		\begin{scope}[xshift=-0.5cm,yshift=-7cm]
		\fill[gray!20] (2,0) rectangle (3,5);
		\fill[gray!20] (0,2) rectangle (5,3);
		\draw  (0,0) grid +(5,5);
		\draw[thick] 
		(0,0) -- (5,0) -- (5,5) -- (0,5) -- (0,0);
		\filldraw[brown!50!black]
		(0.5, 4.5) circle (8pt)
		(1.5, 3.5) circle (8pt)
		(2.5, 2.5) circle (8pt)
		(3.5, 1.5) circle (8pt)
		(4.5, 0.5) circle (8pt);
		\draw[thin, dashed, color=red!50] (0,0)--(5,5) 
		(5,0) -- (0,5);
		\end{scope}					    
		\begin{scope}[xshift=7.5cm,yshift=-7cm]
		\fill[gray!20] (2,0) rectangle (3,5);
		\fill[gray!20] (0,2) rectangle (5,3);
		\draw  (0,0) grid +(5,5);
		\draw[thick] 
		(0,0) -- (5,0) -- (5,5) -- (0,5) -- (0,0);
		\filldraw[brown!50!black]
		(0.5, 3.5) circle (8pt)
		(1.5, 4.5) circle (8pt)
		(2.5, 2.5) circle (8pt)
		(3.5, 0.5) circle (8pt)
		(4.5, 1.5) circle (8pt);
		\draw[thin, dashed, color=red!50] (0,0)--(5,5) 
		(5,0) -- (0,5);
		\end{scope}					    				    
		\begin{scope}[xshift=15.5cm,yshift=-7cm]
		\fill[gray!20] (2,0) rectangle (3,5);
		\fill[gray!20] (0,2) rectangle (5,3);	
		\draw  (0,0) grid +(5,5);
		\draw[thick] 
		(0,0) -- (5,0) -- (5,5) -- (0,5) -- (0,0);
		\filldraw[brown!50!black]
		(0.5, 4.5) circle (8pt)
		(1.5, 1.5) circle (8pt)
		(2.5, 2.5) circle (8pt)
		(3.5, 3.5) circle (8pt)
		(4.5, 0.5) circle (8pt);
		\draw[thin, dashed, color=red!50] (0,0)--(5,5) 
		(5,0) -- (0,5);
		
		\end{scope}				
		
		\end{tikzpicture}
		\caption{The three $4\x 4$ and $5\x 5$ inequivalent diagonally symmetric rook placements.}
		\label{boards3}
	\end{center}
\end{figure}
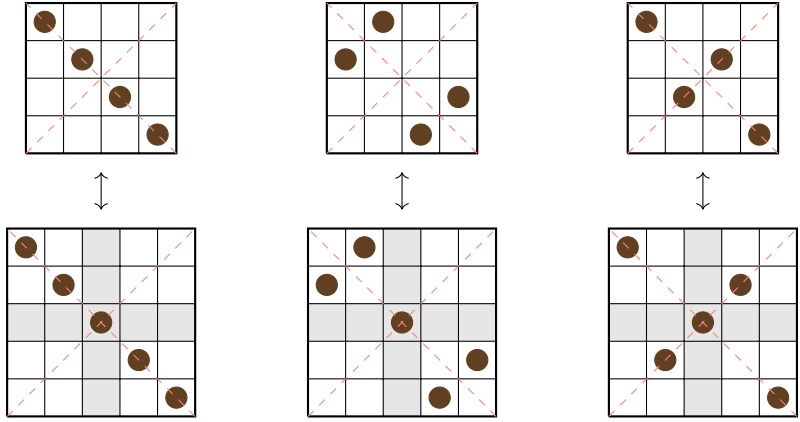
From now on, let $d_n$ denote the number of inequivalent diagonally symmetric arrangements of $n$ non-attacking rooks on an $n \x n$ board. 
\begin{nlem}
	For all $n\geq 1$, $d_{2n}=d_{2n+1}$.
\end{nlem}
\begin{proof}
	Note that the center square must contain a rook in the odd case; the result follows by deleting the middle row and column from a $(2n + 1) \x (2n + 1)$ rook placement.
\end{proof}
Figure \ref{boards3} illustrates an instance of this correspondence, where $n=2$. 

Thus, one only really needs to count diagonally symmetric placements on even side-length boards. We will prove that $\del_n=d_{2n}$ by exhibiting an explicit bijection between $DIII$ clans and diagonally symmetric rook placements. 

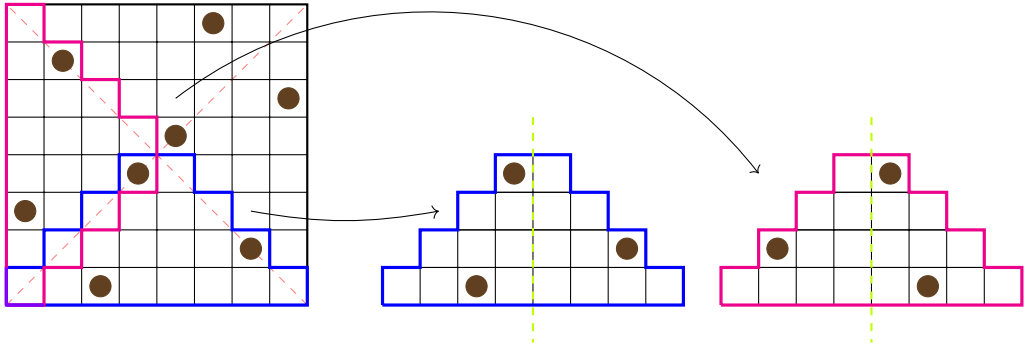
\begin{figure}[!ht] 
	\begin{center}
		\begin{tikzpicture}[scale=0.5]
		\draw[thin, dashed, color=red!50] (0,0)--(8,8)  (8,0) -- (0,8);
		\draw  (0,0) grid +(8,8);
		\draw[thick] 
		(0,0) -- (8,0) -- (8,8) -- (0,8) -- (0,0);
		\draw[very thick, color=blue] (0,0) -- (0,1) -- (1,1) -- (1,2) -- (2,2) -- (2,3) -- (3,3) -- (3,4) -- (5,4) -- (5,3) -- (6,3) -- (6,2) -- (7,2) -- (7,1) -- (8,1) -- (8,0) -- (0,0);
		\draw[very thick, color=magenta] (0,0) -- (1,0) -- (1,1) -- (2,1) -- (2,2) -- (3,2) -- (3,3) -- (4,3) -- (4,5) -- (3,5) -- (3,6) -- (2,6) -- (2,7) -- (1,7) -- (1,8) -- (0,8) -- (0,0);
		\draw[very thick, color=blue!50!magenta] (0,1) -- (0,0) -- (1,0);
		\filldraw[brown!50!black] (0.5, 2.5) circle (8pt)
		(1.5, 6.5) circle (8pt)
		(2.5, 0.5) circle (8pt)
		(3.5, 3.5) circle (8pt)
		(4.5, 4.5) circle (8pt)
		(5.5, 7.5) circle (8pt)
		(6.5, 1.5) circle (8pt)
		(7.5, 5.5) circle (8pt);
		\draw[->] (6.5,2.5) to [bend right=10] (11.5,2.5);
		\draw[->] (4.5,5.5) to [bend left=45] (20,3.5);
		\begin{scope}[xshift=14cm]
		\newcount\pyramidRow   
		\foreach \Row in {4,...,1} {
			\draw (0,\the\pyramidRow) +(-\Row,0) grid ++(\Row,1);
			\global\advance\pyramidRow by 1  
		}
		\draw[very thick, color=blue, xshift=-4cm] (0,0) -- (0,1) -- (1,1) -- (1,2) -- (2,2) -- (2,3) -- (3,3) -- (3,4) -- (5,4) -- (5,3) -- (6,3) -- (6,2) -- (7,2) -- (7,1) -- (8,1) -- (8,0) -- (0,0);
		\filldraw[brown!50!black] 
		(-0.5,3.5) circle (8pt)
		(-1.5,0.5) circle (8pt)
		(2.5,1.5) circle (8pt);
		\draw[lime, dashed, thick] (0,5) -- (0,-1);	
		\end{scope}
		\begin{scope}[xshift=23cm]
		\newcount\pyramidRow   
		\foreach \Row in {4,...,1} {
			\draw (0,\the\pyramidRow) +(-\Row,0) grid ++(\Row,1);
			\global\advance\pyramidRow by 1  
		}
		\draw[very thick, color=magenta, xshift=-4cm] (0,0) -- (0,1) -- (1,1) -- (1,2) -- (2,2) -- (2,3) -- (3,3) -- (3,4) -- (5,4) -- (5,3) -- (6,3) -- (6,2) -- (7,2) -- (7,1) -- (8,1) -- (8,0) -- (0,0);
		\filldraw[brown!50!black] 
		(0.5,3.5) circle (8pt)
		(1.5,0.5) circle (8pt)
		(-2.5,1.5) circle (8pt);
		\draw[lime, dashed, thick] (0,5) -- (0,-1);	
		\end{scope}
		\end{tikzpicture}
		\caption{The two possible pyramids of a diagonally symmetric rook placement.}
		\label{pyramid1}
	\end{center}
\end{figure}

The diagonals of a $2n \x 2n$ board divide it into four triangles, and one sees that the information of a diagonally symmetric rook placement is captured within any of these triangles.\footnote{A diagonally symmetric rook placement (equivalence class) produces two possible pyramids which differ by a reflection. See Figure~\ref{pyramid1}.} In identifying rook placements, we will then extract one of the triangles of the board (rooks included), to obtain a \emph{pyramid}. It will be convenient to introduce coordinates on the blocks of the pyramids by dividing them into \emph{left} and \emph{right} halves. The indices on the left $l_{i,j}$ increase moving up and to the right, while those on the right $r_{i,j}$ increase as we move up and to the left (see Figure~\ref{pyramid2}). The following characterization is self-evident.

\begin{nlem} \label{pyramidChar} A pyramid corresponds to a diagonally symmetric rook placement if and only if for each \hbox{$1\leq k \leq n$}, there is a unique block of the pyramid, $l_{i,j}$ or $r_{i,j}$, on which a rook is placed and for which either $i=k$, or $j=k$, or both $i=j=k$.  
\end{nlem}
%\begin{proof}
%Clearly there can be at most one rook in each row and column of the pyramid. Starting with a rook placement with rook at $b_{i,j}$, reflect it around the board and further rule out the rows and columns of the resulting rooks. Tracing these prohibited rows and columns into the pyramid gives the statement in terms left and right pyramid coordinates.
%\end{proof}

Now we give an algorithm for obtaining a pyramid from a clan by reading the symbols $c_1$ through $c_n$ in reverse order, and placing rooks as we descend rows of the pyramid. An auxiliary variable $X$ acts as a ``switch'' between the left and right sides of the pyramid; every time we encounter a $-$ or the first mate of a $\Pi_1$ pair, the switch gets ``flipped.'' 

\begin{nalgom} \label{pyrAlgo} Given a $DIII$ $(n,n)$-clan $\gamma$, we construct a pyramid corresponding to a diagonally symmetric $2n \x 2n$ rook placement equivalence class as follows.

set $i=n$, $X=l$. 

while $i \geq 1$:

\qquad if $c_i =+$:

\qquad \qquad place rook at $X_{i,i}$

\qquad if $c_i = -$:

\qquad \qquad flip switch

\qquad \qquad place rook at $X_{i,i}$

\qquad if $c_i \in \N$:

\qquad \qquad find $c_j=c_i$ in $\gamma$

\qquad \qquad if $j>n$ and $2n+1-j>i$: \hspace{1.5cm} [second condition prevents redundacy]

\qquad \qquad \qquad place a rook at $X_{i,2n+1-j}$ \hspace{1.5cm} [indices corr. to opposing $\Pi_0$ pairs]

\qquad \qquad if $i<j\leq n$: \hspace{5.5cm} [positions of a $\Pi_1$ pair] 

\qquad \qquad \qquad flip switch 

\qquad \qquad \qquad place rook at $X_{i,j}$ 

\qquad \qquad else:

\qquad \qquad \qquad pass 

\qquad subtract 1 from i
	
\end{nalgom}
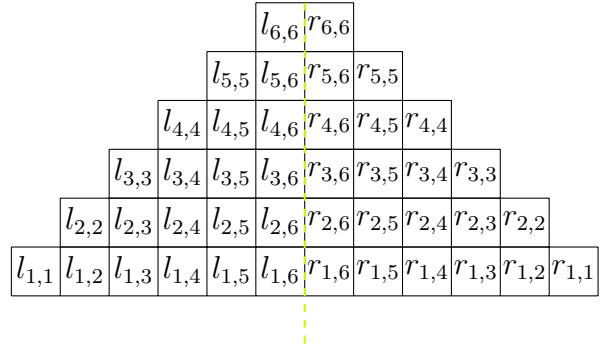
\begin{figure}[tp!] 
	\caption{Coordinates on pyramids.} \label{pyramid2}
	\begin{center}
		\begin{tikzpicture}[scale=0.65]
		\pyramid{6}
		\foreach \x in {1,...,6} { 
			\foreach \y in {\x,...,6} {
				\draw (\y-6.5, \x-0.5) node   {$l_{\x,\y}$};
			}
		}
		\foreach \x in {1,...,6} { 
			\foreach \y in {6,...,\x} {
				\draw (-\y+6.5, \x-0.5) node   {$r_{\x,\y}$};
			}
		}
		\draw[lime, dashed, thick] (0,6) -- (0,-1);		
		\begin{comment}
		\draw[<->] (5,3) to (9,3);
		
		\begin{scope}[xshift=14cm]
		\pyramid{6}
		\foreach \x in {1,...,6} { 
			\foreach \y in {\x,...,6} {
				\draw (\y-6.5, \x-0.5) node   {$b_{\x,\y}$};
			}
		}
		\foreach \x in {6,...,1} { \pgfmathsetmacro\xlab{7-\x}
			\foreach \y in {1,...,\x} {
				\pgfmathsetmacro\ylab{\y+6}
				\draw (\y-0.5, 6.5-\x) node [font=\small]  {$b_{\pgfmathprintnumber{\xlab},\pgfmathprintnumber{\ylab}}$};
			}
		}
		\draw[lime, dashed, thick] (0,6) -- (0,-1);		}
		\end{scope}
		\end{comment}
		\end{tikzpicture}
	\end{center}
\end{figure} 

It is easy to verify that this algorithm produces a pyramid that satisfies the condition of Lemma \ref{pyramidChar}, yielding a diagonally symmetric rook placement. As an example, the blue pyramid in Figure~\ref{pyramid1} is obtained from the $(4,4)$-clan, $1{-}1{+}{-}2{+}2$.

Without trouble, this algorithm can be reversed to give a map from pyramids to clans. However, exactly one of the two pyramids from a given rook placement produces a $DIII$ clan. For example, the pink pyramid in Figure~\ref{pyramid1} would yield the clan $1{-}1{-}{+}2{+}2$, which violates condition 3 of Definition~\ref{defn:Dclans}. In general, if one pyramid produces $\gamma$, then the other produces $\Flip(\gamma)$. So each diagonally symmetric rook placement contains a unique pyramid which gives a $DIII$ clan, completing the bijection. 

\begin{nthm} \label{pyrBij}Diagonally symmetric rook placements on a $2n \x 2n$ board and $(n,n)$-clans of type $DIII$ are in bijection, whence $\del_n=d_{2n}$.
\end{nthm}

\begin{nrk}
Consider a $2n \x 2n$ rook placement as a permutation matrix $v$ once again. Symmetry across $d$ implies that $v$ is the matrix of an involution, while symmetry across $d'$ implies that $v$ is a signed permutation via the usual embedding into $\mc{S}_{2n}$. Then, in the notation of Proposition~\ref{prop:d8}, $R\cdot v = w_0 v$, which is the involution which takes $i\mapsto 2n+1-v(i)$. In terms of signed permutations of $[n]$, $(R\cdot v) (i) = -v(i)$.\footnote{This also implies that $R^2 \cdot v = w_0vw_0=v$.} Thus, diagonally symmetric rook placements up to equivalence (and $DIII$ clans) are also in bijection with pairs of signed permutations $\{v, R\cdot v\}$ of order two. We thank one of the anonymous referees for pointing this out to us.
\end{nrk}

\begin{comment}
We will call a value $i$ such that $v(i)=-i$ (or $v(i)=2n+1-i$, as an element of $\mc{S}_{2n}$) a \emph{negation point} of $v$.

\begin{proof}
In view of the preceding discussion, it suffices to exhibit a bijection between $DIII$ clans and pairs of involutions $\{v, R\cdot v\}$ which are signed permutations. Between $v$ and $R\cdot v$, exactly one has the following property: the number of transpositions of the form $(i, j)$ where $1\leq i,j \leq n$ plus the number of negation points is even.

\end{comment}
\subsection{Minimally Intersecting Set Partitions} \label{subsec:setparts}
Consider partitions of the set $[n]=\{1,\dots,n\}$ ordered by refinement. Two partitions $p$ and $p'$ are said to be \emph{minimally intersecting} if the partition whose blocks are the pairwise intersections of blocks from $p$ and $p'$ is the minimal partition
$$\mbox{$p_{\min}= \{ \{1\}, \{2\}, \dots ,\{n\} \}$}.$$ Lemma 2 of~\cite{pittel} says that the right hand side of the equation (\ref{eq:dgen}) is equal to $e^x$ plus the exponential generating function for the number of ordered pairs of minimally intersecting partitions $(p,p')$ of $[n]$ such that $p$ consists of exactly two blocks. In other words, the number of $DIII$ $(n,n)$-clans is one more than the number of such pairs of partitions \cite[\href{https://oeis.org/A000902}{A000902}]{OEIS}. In this subsection, we present a map between pyramids and partition pairs that demonstrates this equality. 

\begin{nrk} There is a well known bijection between $n \x n$ staircase rook placements and partitions of $[n+1]$ (see~\cite[pp. 77-78]{loehr}). Observing that each pyramid consists of two staircase shapes with ``complementary'' rook placements, one could also define a correspondence with pairs of partitions of $[n+1]$ with certain properties. We leave this description to the motivated reader. 
\end{nrk}

Consider a pyramid that corresponds to an $(n,n)$-clan which is not ${+} \cdots {+}{-}\cdots {-}$. First we describe how to obtain the two-block partition of the corresponding pair, which we will write as $p=\{L, R\}$. 
\begin{enumerate}
\item If there is a rook at $l_{i,i}$, then $i\in L$; if there is a rook at $r_{i,i}$, then $i\in R$. 
\item If there is a rook at $l_{i,j}$ for $i\neq j$, then $j\in L$ and $i \in R$. Similarly, if there is a rook at $r_{i,j}$, then $j\in R$ and $i\in L$.
\end{enumerate}
Then we construct $p'$ by taking $\{i,j\}$ as a block for each rook at $l_{i,j}$ or $r_{i,j}$. Thus, the blocks of $p'$ have maximum size two, and rooks at $l_{i,i}$ or $r_{i,i}$ give blocks that are singletons. It is clear that the pair $(p,p')$ is minimally intersecting.

\begin{nexap} The blue pyramid of Figure~\ref{pyramid1} gives the pair $(p,p')$ with $p=\{\{3,4\}, \{1,2\}\}$ and $p'=\{ \{1,3\}, \{2\}, \{4\}\}$.
\end{nexap}

Notice that reflecting a pyramid across the center line swaps the blocks $L$ and $R$ of $p$, but $(p,p')$ is unchanged. Exclusion of the clan ${+}\cdots{+}{-}\cdots{-}$ guarantees that neither $L$ nor $R$ is empty.

%Given a pyramid, to guarantee that neither of $L$ or $R$ is empty, we only need one rook placed at some $l_{i,j}$ or $r_{i,j}$ where $i \neq j$, or if all of the rooks are at $l_{i,i}$'s or $r_{i,i}$'s (which occurs when the underlying clan consists only of $+$'s and $-$'s), then at least one rook needs to be on each half. This only fails for the clan $+\cdots +-\cdots -$, which gives a pyramid with rooks only at $l_{i,i}$ for all $1\leq i \leq n$. In all other cases, the rook placements on the pyramids clearly determine unique pairs of minimally intersecting partitions. 

Now we describe how to obtain a pyramid from a pair $(p,p')$, where $p=\{L,R\}$. Observe that $p'$ cannot have any blocks of size greater than two.%; suppose otherwise that $b$ was a block in $p'$ of size three or greater. Then, by the pigeonhole principle we would have at least two of the elements of $b$ in either $L$ or $R$. The pairwise intersections of these blocks would leave us with a block of size at least two, and so $p$ and $p'$ could not be minimally intersecting. With this in mind, we construct a pyramid from a pair $(p,p')$ as follows.
\begin{enumerate}
\item If $i\in L$ (respectively, in $R$) and $i$ is a singleton in $p'$, then place a rook at $l_{i,i}$ (respectively, at $r_{i,i}$.%; if $i\in R$ and $i$ is a singleton in $p'$, then place a rook at $r_{i,i}$.
\item If $j \in L$ (respectively, in $R$) and $\{i,j\}$ is a block of $p'$ with $i<j$, then place a rook at $l_{i,j}$ (respectively, at $r_{i,j}$).%; if $j \in R$ and $\{i,j\}$ is a block of $p'$ with $i<j$ then place a rook at $r_{i,j}$.
\end{enumerate}

%It is clear that the pyramid constructed will satisfy the criterion of Lemma \ref{pyramidChar}, and so it corresponds to a unique diagonally symmetric rook placement and type $DIII$ clan. 
This recipe inverts the (partial) map from pyramids to partition pairs described above, establishing the following. %. Again, the pyramid corresponding to the clan $+\cdots +-\cdots -$ can not be constructed through this process. This discussion establishes the following.

\begin{nthm} The set of $DIII$ $(n,n)$-clans without the clan ${+}\cdots{+}{-}\cdots{-}$ is in bijection with the set of ordered pairs $(p,p')$ of minimally intersecting pairs of partitions of $[n]$, where $p$ has exactly two blocks.
\end{nthm}

\subsection{Lattice paths}
\label{subsec:LP}
%In this final subsection, we present a construction similar to those made in~\cite{canGenesis} and~\cite{aoSects} for type $AIII$ and type $CI$, respectively. %While $DIII$ clans can be viewed as subsets of the clans of either of those types, the weighted lattice paths presented here are only a subset of those presented in the latter work. 

Recall that an \emph{$(n,n)$ Delannoy path} is an integer lattice path from $(0,0)$ to $(n,n)$ in the plane $\R^2$ consisting only of single north, east, or diagonally northeast steps. Alternatively, one can consider strings from the alphabet $\{N, E, D \}$ such that the number of $N$'s plus the number of $D$'s is equal to the sum of the numbers of $E$'s and $D$'s (which is equal to $n$). We will demonstrate a bijection between the set of $DIII$ \mbox{$(n, n)$-clans} and the set of $(n,n)$ Delannoy paths with certain labels which are defined as follows.

\begin{ndefn} \label{defn:paths}
	By a \emph{labelled step} we mean a pair $(L,l)$,
	where $L\in \{N,E,D\}$ and $l$ is a positive integer
	such that $l=1$ if $L=N$ or $L=E$. 
	A \emph{weighted $(n,n)$ Delannoy path} is a word of the form 
	$W:=W_1\dots W_r$, where the $W_i$'s are labeled steps $W_i=(L_i,l_i)$ such that 
	\begin{enumerate}
		\item $L_1\dots L_r$ is an $(n,n)$ Delannoy path;
		\item $L_i=N$ if and only if $L_{r+1-i}=E$;
		\item letting $k_i=\# \{j < i \mid l_j\neq 1\}$, if $L_i=D$ then $2\leq l_i\leq 2n+1-2(i+2k_i)$ for $1\leq i \leq \floor{\frac{r}{2}}$, and $W_{r+1-i}=(D,2n+3-2(i+2k_i)-l_i)$;
		\item either $L_\frac{r}{2}=E$ (so that $L_{\frac{r}{2}+1}=N$) or $W_\frac{r}{2}=(D, 3)$ (so that $W_{\frac{r}{2}+1}=(D,2)$).\footnote{As a consequence of the preceding properties, $r$ is guaranteed to be even.}
	\end{enumerate}
	\end{ndefn}	

\begin{nthm}
	There is a bijection between the set of weighted $(n, n)$ Delannoy 
	paths and the set of DIII $(n,n)$-clans. 
%	In particular, we have 
%	$$
%	\Delta_{n} = \sum_{W\in \mc{D}^{\omega}(n)} 1.
%	$$
\end{nthm}

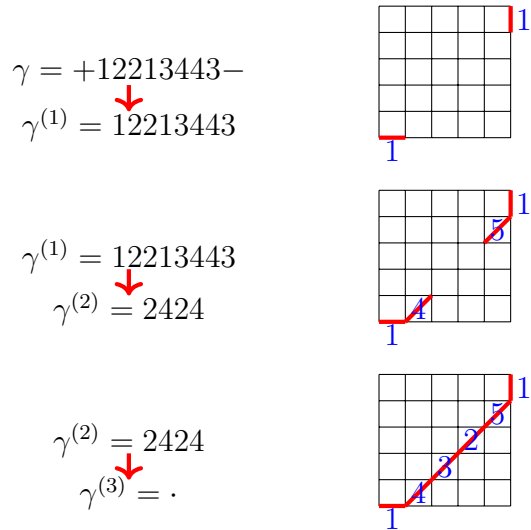
\begin{figure}[htbp!]\label{5:lastpic}
		\begin{center}
			\begin{tikzpicture}[scale=.35]
			\begin{scope}[xshift = -5cm]
			\node at (0,2.5) {$\gamma =+12213443-$};
			\node at (0,.5) {$\gamma^{(1)}= 12213443$};
			\draw[ultra thick,red, ->] (0,2) to (0,1);
			\end{scope}
			\begin{scope}[xshift = 4.5cm]
			\draw (0,0) grid (5,5);
			\draw[ultra thick,red] (5,5) to (5,4);
			\draw[ultra thick,red] (0,0) to (1,0);
			\node[blue] at (5.5, 4.5) {$1$};
			\node[blue] at (0.5, -0.5) {$1$};
			\end{scope}

			\begin{scope}[xshift = -5cm,yshift=-7cm]
			\node at (0,2.5) {$\gamma^{(1)}= 12213443$};
			\node at (0,.5) {$\gamma^{(2)}= 2424 $};
			\draw[ultra thick,red, ->] (0,2) to (0,1);
			\end{scope}
			\begin{scope}[xshift = 4.5cm,yshift=-7cm]
			\draw (0,0) grid (5,5);
			\draw[ultra thick,red] (5,5) to (5,4);
			\draw[ultra thick,red] (5,4) to (4,3);
			\draw[ultra thick,red] (0,0) to (1,0);
			\draw[ultra thick,red] (1,0) to (2,1);
			\node[blue] at (5.5, 4.5) {$1$};
			\node[blue] at (0.5, -0.5) {$1$};		
			\node[blue] at (4.5, 3.5) {$5$};
			\node[blue] at (1.5, .5) {$4$};
			\end{scope}

			\begin{scope}[xshift = -5cm,yshift=-14cm]
			\node at (0,2.5) {$ \gamma^{(2)}= 2424$};
			\node at (0,.5) {$\gamma^{(3)}= \cdot$};
			\draw[ultra thick,red, ->] (0,2) to (0,1);
			\end{scope}
			\begin{scope}[xshift = 4.5cm,yshift=-14cm]
			\draw (0,0) grid (5,5);
		\draw[ultra thick,red] (5,5) to (5,4);
		\draw[ultra thick,red] (5,4) to (4,3);
		\draw[ultra thick,red] (0,0) to (1,0);
		\draw[ultra thick,red] (1,0) to (2,1);
		\draw[ultra thick,red] (2,1) to (3,2);
		\draw[ultra thick,red] (3,2) to (4,3);
		\node[blue] at (5.5, 4.5) {$1$};
		\node[blue] at (0.5, -0.5) {$1$};		
		\node[blue] at (4.5, 3.5) {$5$};
		\node[blue] at (1.5, .5) {$4$};
			\node[blue] at (2.5, 1.5) {$3$};
			\node[blue] at (3.5, 2.5) {$2$};
			\end{scope}
			\end{tikzpicture}
		\end{center}
		\caption{Algorithmic construction of the bijection onto weighted Delannoy paths.}
		\label{5:last pic}
	\end{figure}
\begin{proof} We will indicate how to obtain a weighted $(n,n)$ Delannoy path from a type $DIII$ clan $\gamma=c_1 \cdots c_{2n}$. If $c_{2n}$ is a $-$ sign, we draw an $N$-step from $(n, n-1)$ to $(n, n)$ and an $E$-step between $(0, 0)$ and $(1, 0)$. Then we remove $c_1$ and $c_{2n}$ from $\gamma$ to obtain $\gamma^{(1)}=c_2\cdots c_{2n-1}$.

In a similar manner, if $c_{2n}=+$, we draw an $E$-step from \mbox{$(n-1, n)$} to $(n, n)$ and an $N$-step between $(0, 0)$ and $(0, 1)$. Again we remove $c_1$ and $c_{2n}$ from $\gamma$, but in this case we then swap $c_n$ and $c_{n+1}$ to obtain $\gamma^{(1)}=c_2 \cdots c_{n+1}c_n\cdots c_{2n-1}$.
	
	If $c_{2n}$ is a natural number from $\Pi_0$ pair $(c_i, c_{2n})$ ($i\leq n$), we draw a $D$-step between $(n-1, n-1)$ and $(n, n)$ and label this step $i$, and we draw another $D$-step between $(0, 0)$ and $(1,1)$ and label this step $2n+1-i$. Then we remove all four symbols $c_1, c_i, c_{2n+1-i}$, and $c_{2n}$ from $\gamma$ and call the resulting $(n-2,n-2)$-clan $\gamma^{(1)}$.
	
	In case $c_{2n}$ is a natural number from $\Pi_1$ pair $(c_j,c_{2n})$ ($j>n$) with opposing $\Pi_1$ pair $(c_1,c_{2n+1-j})$, then we draw a $D$-step between $(n-1, n-1)$ and $(n, n)$ and label this step $j$, and we draw another $D$-step between $(0, 0)$ and $(1,1)$ and label this step with $2n+1-j$. We remove all four symbols $c_1, c_j, c_{2n+1-j}$, and $c_{2n}$ from $\gamma$, then swap $c_n$ and $c_{n+1}$ and call the resulting $(n-2,n-2)$-clan $\gamma^{(1)}$.

After performing this first step, we iterate the same procedure upon $\gamma^{(1)}$ by examining its last symbol, thereby obtaining $\gamma^{(2)}$ and so on, building the path from the corners inwards. 

This is clearly an injective construction. The complicated condition 3 of Definition~\ref{defn:paths} on the weights (which give the placement of mates in the clan) just guarantees that the construction can be reversed to obtain a skew-symmetric clan. Condition 4 guarantees the parity condition of Definition~\ref{defn:Dclans}, so we have a bijection. 
\end{proof}

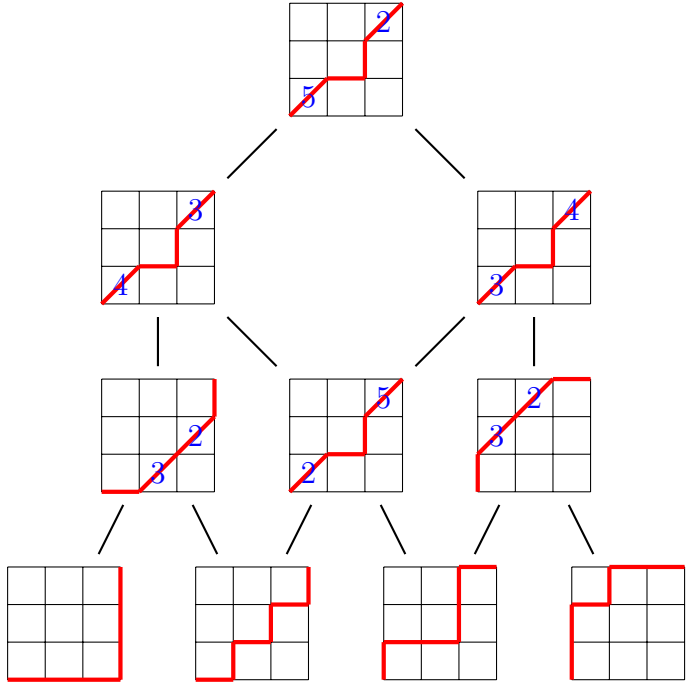
\begin{figure}[h!]
	\begin{center}
		\begin{tikzpicture}[scale=.25]
		\node at (-15,0) (a1) {$\begin{tikzpicture}[scale=.5]
			\draw (0,0) grid (3,3);	
			\draw[ultra thick,red] (3,3) to (3,0);
			\draw[ultra thick,red] (0,0) to (3,0);
	\end{tikzpicture}	$};
		\node at (-5,0) (a2) {$\begin{tikzpicture}[scale=.5]
			\draw (0,0) grid (3,3);
			\draw[ultra thick,red] (0,0) to (1,0);
			\draw[ultra thick,red] (1,0) to (1,1);
			\draw[ultra thick,red] (1,1) to (2,1);
			\draw[ultra thick,red] (2,1) to (2,2);
			\draw[ultra thick,red] (2,2) to (3,2);
			\draw[ultra thick,red] (3,2) to (3,3);
			\end{tikzpicture}	$};
		\node at (5,0) (a3) {$\begin{tikzpicture}[scale=.5]
			\draw (0,0) grid (3,3);
			\draw[ultra thick,red] (0,0) to (0,1);
\draw[ultra thick,red] (0,1) to (2,1);
\draw[ultra thick,red] (2,1) to (2,3);
\draw[ultra thick,red] (2,3) to (3,3);
			\end{tikzpicture}	$};
		\node at (15,0) (a4) {$\begin{tikzpicture}[scale=.5]
			\draw (0,0) grid (3,3);
\draw[ultra thick,red] (0,0) to (0,2);
\draw[ultra thick,red] (0,2) to (1,2);
\draw[ultra thick,red] (1,2) to (1,3);
\draw[ultra thick,red] (1,3) to (3,3);
			\end{tikzpicture}	$};

		\node at (-10,10) (b1) {$\begin{tikzpicture}[scale=.5]
			\draw (0,0) grid (3,3);
			\draw[ultra thick,red] (0,0) to (1,0);
			\draw[ultra thick,red] (3,3) to (3,2);
			\draw[ultra thick,red] (3,2) to (1,0);
			\node[blue] at (1.5,.5) {$3$};
			\node[blue] at (2.5,1.5) {$2$};
			\end{tikzpicture}	$};
		\node at (0,10) (b2) {$\begin{tikzpicture}[scale=.5]
			\draw (0,0) grid (3,3);
			\draw[ultra thick,red] (0,0) to (1,1);
			\draw[ultra thick,red] (3,3) to (2,2);
			\draw[ultra thick,red] (2,2) to (2,1);
			\draw[ultra thick,red] (2,1) to (1,1);
			\node[blue] at (.5,.5) {$2$};
			\node[blue] at (2.5,2.5) {$5$};
			\end{tikzpicture}	$};
		\node at (10,10) (b3) {$\begin{tikzpicture}[scale=.5]
			\draw (0,0) grid (3,3);
			\draw[ultra thick,red] (0,0) to (0,1);
\draw[ultra thick,red] (3,3) to (2,3);
\draw[ultra thick,red] (2,3) to (0,1);
\node[blue] at (.5,1.5) {$3$};
\node[blue] at (1.5,2.5) {$2$};
			\end{tikzpicture}	$};

		\node at (-10,20) (c1) {$\begin{tikzpicture}[scale=.5]
			\draw (0,0) grid (3,3);
		\draw[ultra thick,red] (0,0) to (1,1);
\draw[ultra thick,red] (3,3) to (2,2);
\draw[ultra thick,red] (2,2) to (2,1);
\draw[ultra thick,red] (2,1) to (1,1);
\node[blue] at (0.5,0.5) {$4$};
\node[blue] at (2.5,2.5) {$3$};
			\end{tikzpicture}	$};

		\node at  (10,20) (c2) {$\begin{tikzpicture}[scale=.5]
			\draw (0,0) grid (3,3);
		\draw[ultra thick,red] (0,0) to (1,1);
		\draw[ultra thick,red] (3,3) to (2,2);
		\draw[ultra thick,red] (2,2) to (2,1);
		\draw[ultra thick,red] (2,1) to (1,1);
		\node[blue] at (0.5,0.5) {$3$};
		\node[blue] at (2.5,2.5) {$4$};
			\end{tikzpicture}	$};
		
		\node at (0,30) (d) {$\begin{tikzpicture}[scale=.5]
			\draw (0,0) grid (3,3);
			\draw[ultra thick,red] (0,0) to (1,1);
			\draw[ultra thick,red] (3,3) to (2,2);
			\draw[ultra thick,red] (2,2) to (2,1);
			\draw[ultra thick,red] (2,1) to (1,1);
			\node[blue] at (0.5,0.5) {$5$};
			\node[blue] at (2.5,2.5) {$2$};
			\end{tikzpicture}	$};

		\draw[-,  thick] (a1) to (b1);
		\draw[-,  thick] (a2) to (b1);
		\draw[-,  thick] (a2) to (b2);
		\draw[-,  thick] (a3) to (b2);
		\draw[-,  thick] (a3) to (b3);
		\draw[-,  thick] (a4) to (b3);
		
		\draw[-, thick] (b1) to (c1);
		\draw[-, thick] (b2) to (c1);
		\draw[-, thick] (b2) to (c2);
		\draw[-, thick] (b3) to (c2);

		\draw[-, thick] (c1) to (d);
		\draw[-, thick] (c2) to (d);
		
		\end{tikzpicture}
		
		\caption{Weak order on weighted $(3,3)$ Delannoy paths}
		\label{5:LPaths}
	\end{center}
\end{figure}

\begin{nexap}
	Let $\gamma = {+} 1 2 2 1 3 4 4 3 {-}$.
	The steps of our construction 
	are shown in Figure~\ref{5:last pic}.
	\end{nexap}
To supply further examples, we depict the weighted Delannoy paths corresponding to $DIII$ (3,3)-clans in Figure~\ref{5:LPaths} in their weak order poset. 

\section{The big sect} \label{S:bigsect}

In this section, we investigate the number of $(n,n)$-clans in the largest sect; we denote this number by $\epsilon_n$.  These are the clans whose corresponding $B$-orbits comprise the preimage of the dense Schubert cell under the map $\pi:G/L \to G/P$. Since this sect must include the dense $B$-orbit corresponding to the clan $\gamma_0$ of (\ref{eq:gamma0}), we see that this sect has base clan $$\underbrace{{-} {-} \cdots  {-} {-}}_{\text{first}\ n\text{-spots}} {+} {+}\cdots  {+} {+} \;\; \text{or} \;\; \underbrace{{-} {-} \cdots  {-} {+}}_{\text{first}\ n\text{-spots}} {-} {+}\cdots {+}{+} \;\;  $$
depending on whether $n$ is even or odd, respectively. Consequently, a clan lies in the largest sect only if 
\begin{enumerate}[label=(\alph*)]
	\item it has natural number pairs only in $\Pi_0$ when $n$ is even, 
	\item or it has at most two $\Pi_1$ pairs at $(c_{i},c_n)$ and $(c_{2n+1-i},c_{n+1})$ when $n$ is odd.
\end{enumerate}

If $2r$ is the number of pairs of matching natural numbers in a $DIII$ clan $\gamma = c_1 \cdots c_{2n}$ which lies in the largest sect, the clan is determined by pairing $2r$ of the symbols among $c_1\cdots c_{n}$. This can be done in ${n\choose 2r}\frac{(2r)!}{r!2^r}$ many different ways. %Then, the $2r$ positions form $r$ pairs\footnote{These paired positions $(i,j)$ determine first symbols of opposing $\Pi_0$ pairs unless $n$ is odd and $j=n$ in which case the pairing forces $(c_i,c_j)\in \Pi_1$.}, which can be done in $\frac{(2r)!}{2^{r} r!}$ different ways. 
Summing over possible values for $r$, we have
\begin{align}\label{sectnumber}
\epsilon_n = \sum_{r=0}^{\floor{n/2}} \frac{n!}{(n-2r)! r!2^{r}},
\end{align} 
which happens to be the number of involutions on $n$ letters \cite[\href{https://oeis.org/A000085}{A000085}]{OEIS}. This coincidence reveals the following.
\begin{nprop}
	Taking $\epsilon_0=1$, the number of clans in the largest sect satisfies the recurrence relation
	\begin{equation}
	\epsilon_n =\epsilon_{n-1} + (n-1)\epsilon_{n-2},
	\end{equation}
	and has exponential generating function 
	\begin{equation}
	\sum_{n=0}^\infty \epsilon_n\frac{x^n}{n!} = e^{x+\frac{x^2}{2}}.
	\end{equation}
\end{nprop}

Recall that a \emph{partial permutation} is a
map $x : \{1, \dots, m\} \lra \{0, \dots, n\}$ satisfying: 
\begin{itemize}
	\item if $x(i) = x(j)$ and $x(i) \neq 0$, then $i = j$.
\end{itemize}
A partial permutation $x$ can be represented by an $m\x n$ matrix $(x_{ij})$, where $x_{ij}$ is 1 if and only if $x (i) = j$ and is
0 otherwise. Note that under this convention we view our matrices as acting on vectors from the right. These partial permutations are also sometimes called \emph{rook placements}, and in case $m =n$, they form a monoid under matrix multiplication called the \emph{rook monoid} and denoted $\mc{R}_n$. 

\begin{ndefn}
	A \emph{partial involution} on $n$ elements is a partial permutation which is represented by a symmetric $n \x n$ matrix. A partial involution with no fixed points is called a \emph{partial fixed-point-free involution}, and the set of such partial involutions is denoted $\mc{PF}_n$.
\end{ndefn}
%We refer to~\cite{bagnocongruence} and~\cite[Chapter 15]{CCA} for relevant theory on partial involutions and partial permutations.

There is a bijection between $\mc{PF}_n$ and the set of invertible involutions $\mc{I}_n$ as follows: a partial involution matrix can be completed to the matrix of an involution by placing a 1 on the diagonal of any row/column without a 1. However, we will prove that $\epsilon_n = \ab{\mc{PF}_n}$ by exhibiting an explicit bijection between
the clans in the largest sect and the partial fixed-point-free involutions.

Let $\gamma=c_1\cdots c_{2n}$ lie in the largest sect. Construct the associated $x \in \mc{PF}_n$ as follows.

\begin{enumerate}
	\item[(i)] If $c_i = \pm$, then take $x(i) = 0$ for all $1 \leq i\leq n$.
	\item[(ii)] If $(c_i, c_j)  \in \Pi_0$, then take $x(i) = 2n+1-j$ and $x({2n+1-j}) = i$ for all $1 \leq i \leq n < j \leq 2n$.
	\item[(iii)] If $(c_{i}, c_n)  \in \Pi_1$, then take $x(i) = n$ and $x(n)= i$.
\end{enumerate}

It is easy to show that this map is invertible. Let us start with a partial fixed-point-free involution $x\in \mc{PF}_n$ and determine its associated $(n, n)$-clan by assigning the first $n$ symbols and then completing using skew-symmetry.

\begin{enumerate}
	\item[(i)] If $x(i) = 0$, then take $c_i = -$ unless $i=n$ and $n$ is odd, in which case $c_i=+$.
	\item[(ii)] If $x(i) = j$ and $x(j) = i$ with $ i  < j $, then we take $(c_i, c_{2n+1-j})  \in \Pi_0$ unless $j=n$ and $n$ is odd in which case we take $(c_{i},c_n)\in\Pi_1$.
\end{enumerate}

This completes proof of the following.
\begin{nthm} Partial fixed-point-free involutions on $n$ letters and $DIII$ $(n,n)$-clans in the largest sect are in bijection.
\end{nthm}

The elements of $\mc{PF}_n$ also parameterize the congruence orbits of the invertible upper triangular $n \x n$ matrices on the skew-symmetric $n\x n$ matrices (with complex entries). This endows them with a poset structure which is the containment order of the corresponding orbit closures, studied in~\cite{cherniavskySkew} and~\cite{cantim}. Let this poset be denoted $(\mc{PF}_n, \leq_{con})$. The order relation $\leq_{con}$ admits a simple combinatorial description in terms of rank-control matrices.

In~\cite{wyserKorbit}, it is pointed out that the full closure (Bruhat) order on $DIII$ $(n,n)$-clans fails in general to be the restriction of the Bruhat order on all $(n,n)$-clans. In particular, Wyser points out that $DIII$ clans $1{+}{-}12{+}{-}2$ and $12341234$ are not related in the Bruhat order in type $DIII$ (as can be observed in Figure~\ref{fig:del4}), though they are related as $AIII$ clans.

Nevertheless, as stated in the introduction, the closure order on $DIII$ clans restricted to the big sect does coincide with that of $(\mc{PF}_n, \leq_{con})$ via the bijection given here. The poset $(\mc{PF}_n, \leq_{con})$ is itself the restriction of the closure order on $\mc{R}_n$, which describes the big sect closure poset in type $AIII$ \cite{bcSects}. 

A combinatorial description of the closure order on all $DIII$ $(n,n)$-clans (based on the work \cite{gandiniMaffei}) will also appear in the first author's Ph.D. Thesis. From this description it is apparent that the failure of the clan closure order to restrict from type $AIII$ to type $DIII$ results from the failure of the Bruhat order on $\mc{S}_{2n}$ to restrict to the Bruhat order on the type $D_n$ Coxeter group. But the clan order \emph{does} restrict from type $AIII$ to type $CI$ for similar reasons, answering part of \cite[Conjecture 3.6]{wyserKorbit}. It would be of great interest to determine general geometric conditions which guarantee the restriction of closure orders in similar settings; we leave the reader to consider this question. 

%Despite these obstacles and having obtained similar results in types $AIII$ and $CI$ in~\cite{bcSects} and~\cite{aoSects}, we find it appropriate to close this paper by conjecturing that the Bruhat order on type $DIII$ $(n,n)$-clans, restricted to the big sect, is isomorphic as a poset to $(\mc{PF}_n, \leq_{con})$. 

\vspace{1em}
\textbf{Acknowledgements.}
We thank Mahir Bilen Can and William McGovern for helpful conversations and suggestions. We are indebted to the anonymous referees who made numerous valuable suggestions and criticisms of our paper, greatly improving the exposition.

\bibliographystyle{alpha}

\end{document}